\documentclass[twoside,11pt]{article}

\usepackage{blindtext}

 \usepackage[abbrvbib, preprint]{jmlr2e}

\usepackage{jmlr2e}

\usepackage{amsmath,amsfonts}
\usepackage{bm}
\usepackage{multirow}
\usepackage{threeparttable}

\newcommand{\btheta}{\bm{\theta}}

\newcommand{\mcS}{\mathcal{S}}
\newcommand{\mcF}{\mathcal{F}}
\newcommand{\mcE}{\mathcal{E}}
\newcommand{\mcD}{\mathcal{D}}
\newcommand{\mcG}{\mathcal{G}}

\newcommand{\mcH}{\mathcal{H}}

\newcommand{\mcB}{\mathcal{B}}

\newcommand{\mcZ}{\mathcal{Z}}
\newcommand{\mcT}{\mathcal{T}}

\newcommand{\mbE}{\mathbb{E}}

\newcommand{\mbN}{\mathbb{N}}

\usepackage{lastpage}
\jmlrheading{23}{2022}{1-\pageref{LastPage}}{1/21; Revised 5/22}{9/22}{21-0000}{Author One and Author Two}

\ShortHeadings{SLP}{Zhou and Yang}
\firstpageno{1}

\begin{document}

\title{On a synergistic learning phenomenon in nonparametric domain adaptation}

\author{\name Ling Zhou \email zhouling@swufe.edu.cn \\
       \addr Center of Statistical Research and School of Statistics and Data Science\\
       Southwestern University of Finance and Economics\\
       Chengdu, 610064, China
       \AND
       \name Yuhong Yang \email  yyangsc@tsinghua.edu.cn\\
       \addr Yau Mathematical Science Center\\
       Tshinghua University\\
       Beijing, 100084, China}

\editor{My editor}

\maketitle

\begin{abstract}
Consider nonparametric domain adaptation for regression, which assumes the same conditional distribution of the response given the covariates but different marginal distributions of the covariates.  An important goal is to 
understand how the source data may improve the minimax convergence rate of learning the regression function when the likelihood ratio of the covariate marginal distributions of the target data and the source data are unbounded. A 
previous work of \cite{pathak2022new} show that the minimax transfer learning rate is simply determined by the faster rate of using either the source or the target data alone.  In this paper, we present a new synergistic learning phenomenon (SLP) that the minimax convergence rate based on both data may sometimes be faster (even much faster) than the better rate of convergence based on the source or target data only.  The SLP occurs when and only when the target sample size is smaller (in order) than but not too much smaller than the source sample size in relation to the smoothness of the regression function and the nature of the covariate densities of the source and target distributions. Interestingly,  the SLP happens in two different ways according to the relationship between the two sample sizes. One is that the target data help alleviate the difficulty in estimating the regression function at points where the density of the source data is close to zero and the other is that the source data (with its larger sample size than that of the target data) help the estimation at points where the density of the source data is not small. 
Extensions to handle unknown source and target parameters and smoothness of the regression function are also obtained. 
\end{abstract}


\section{Introduction}
\label{sec:intro}

Consider transfer learning that focuses on the case of covariate shift, also known as domain adaptation, in the regression setting. Suppose we have i.i.d source data $\mcD = \{X_i, Y_i\}_{i=1}^n$ and i.i.d target data 
$\mcD^T = \{X^T_i, Y^T_i\}_{i=1}^{n_T}$ available. In each pair of observation $(X, Y)$, $X \in [0, 1]$ is the covariate and $Y$ is the response. For both source and target data, given the covariate $X$, assume the following data generating scheme:
\[
Y = f(X) + \epsilon, \; \; \; \epsilon \sim N(0, \sigma^2),
\]
where $\sigma$ is  bounded above by $\overline{\sigma} < \infty$ for a known 
$\overline{\sigma} > 0$. Let $f \in \mcH(\beta, \kappa)$, a H$\ddot{\text{o}}$lder class defined on $[0, 1]$, which includes all functions $f$ satisfying $|f^{(l)}(x) - f^{(l)}(y)| \leq \kappa |x - y|^{\beta - l}$ for all $x, y \in [0, 1]$ and $\sup_x|f(x)| \leq \kappa$, where $l = \lfloor \beta \rfloor$, for some positive constants $\beta$ and $\kappa$. 
Under the framework of domain adaptation, we assume the density functions of $X$ in source data and target data are $h(x)$ and $h_T(x)$, respectively, which are assumed to be known. We aim to study the risk of estimating $f$ using the data $\mcD$ and $\mcD^T$ under $h_T(\cdot)$, namely, 
\[
\mbE\int_0^1 \left( \hat{f}_{n, n_T}(x) - f(x)\right)^2 h_T(x) dx, 
\]
where $\hat{f}_{n, n_T}$ is estimated based on data $\mcD$ and $\mcD^T$ together. Clearly, the squared $L_2(h_T)$ error loss above  emphasizes 
that the interest of learning is on the target population as is described by the target density $h_T$. 

In the research area of covariate shift, starting from the theoretical work in \cite{ben2010theory}, a lot of research has focused on deriving upper and/or lower bounds for the target generalization error \citep{cortes2010learning,cortes2011domain,germain2016new,mansour2009domain,redko2020survey}. However, most of work rely largely on refinements of traditional metrics and divergences between distributions, and often 
only yield a coarse view of when transfer is possible or not \citep{kpotufe2021marginal}.  \cite{kpotufe2017lipschitz} considered error bounds under bounded likelihood ratios, namely, $\sup_x \rho(x) := \sup_x \frac{h_T(x)}{h(x)} < C < \infty$  for some constant $C$. \cite{cortes2010learning} derived upper bounds under a type of bounded second moment condition, i.e., $\mbE_{h}\left(\rho^2(x)\right) < C < \infty$ for some constant $C$. Recently, a few papers focus on deriving learning rates when the bounded likelihood ratio or moment conditions do not hold.  
In particular, let $h_T(x)$ be the uniform distribution on the interval $[0, 1]$ and $h(x) = a x^{a - 1}$ for $x \in [0, 1]$ and some $a \geq 2$. Clearly, the likelihood ratio $\sup_{x}\rho(x) = \infty$ and the moment condition fails, i.e., $\mbE_{h}\left[ \left(\frac{h_T(x)}{h(x)}\right)^2 \right] = \infty$ when $a \geq 2$. \cite{kpotufe2021marginal} introduced the notion of a transfer exponent to explicitly characterize the optimal nonparametric rates for classification. Under the uniform and beta target-source density pair, the minimax rate for a classification problem via the computation from the transfer exponent is $\left(n^{\frac{2\beta + 1}{2\beta + a}} + n_T\right)^{-\frac{2\beta}{2\beta + 1}}$. 
 Later, in the regression setting, \cite{pathak2022new} developed a similarity measure to determine the minimax learning rate of $f$ over a class of target-source density pairs. Particularly, for the uniform and beta target-source density pair,  the minimax rate of learning $f$ via the similarity measure is $\left( n^{\frac{2\beta + 1}{2\beta + a - 1}} + n_T\right)^{-\frac{2\beta}{2\beta + 1}}$ derived under the condition $\beta \leq 1$.
They pointed out that the reason why they obtained a sharper rate of minimax learning when $n^{\frac{2\beta +1}{2\beta + a}}$ is larger (in order) than $n_T$
is that the function class of the target-source density pairs induced by the similarity measure is smaller than that induced by the transfer exponent in \cite{kpotufe2021marginal}. From the results of \cite{pathak2022new}, it can be seen that the minimax learning rate is actually determined by the faster rate of using source data or target data only.  This interesting finding prompts two questions for further understanding and insight on domain adaptation in regression: 
\begin{enumerate}
\item They pursue the worst case rate of convergence when $(h, h_T)$ is over a class of distribution pairs that share the same similarity measure. In their derivation of the minimax lower bound, the ``least favorable'' pair $(h, h_T)$ is chosen to be discrete, which is not applicable to 
the situation that both $h$ and $h_T$ are continuous. The first question is: Given a pair of continuous distribution 
$(h, h_T)$, what is the minimax optimal rate of convergence for domain adaptation? Is it still the same as that given for the 
collection of $(h, h_T)$ pairs in \cite{pathak2022new}?

\item The minimax rate $\left( n^{\frac{2\beta + 1}{2\beta + a - 1}} + n_T\right)^{-\frac{2\beta}{2\beta + 1}}$ identified by \cite{pathak2022new} is in some sense disappointing because  the combined data cannot do better than their individual best (i.e., one does equally well by considering only one of the two data sets)  in terms of rate of convergence. The second question is: With both source and target data available, is it possible to integrate the data to build an estimator that converges faster in the minimax sense than the better rate offered by the source and target data separately?

\end{enumerate}

\cite{schmidt2022local} developed a tool to understand the local convergence behavior of local averaging estimation and used it to construct an estimator for covariate shift. They found that for $f \in \mcH(1, \kappa)$ and $a > 2.5$, with the use of target data and source data together, their upper bound of convergence rate in probability gets improved compared to that of using the source data only. After a simple calculation, one sees that it is actually faster than that of using the target data only as well. It remains unknown what the minimax rate of convergence is for the domain adaptation problem for a given pair of $(h, h_T)$.

In this work, for general $n$ and $n_T$ with $\max(n, n_T) \to \infty$, we explicitly give the minimax rate of learning $f$ for both bounded and unbounded likelihood ratios or second order moment conditions. More importantly, we show that when $n_T$ is smaller than but not too much smaller than $n$ in relation to the smoothness $\beta$ and the parameter $a$ that dictates the likelihood ratio between $h_T$ and $h$, there is a synergistic learning phenomenon (SLP), i.e., the optimal learning rate based on the source data and the target data together is faster than that based on the source data or the target data only. For the purpose of highlighting the new SLP with less secondary technicality, in this work, we start with the one-dimensional H$\ddot{\text{o}}$lder class for the regression function, uniform target density and a general beta source density, and Gaussian errors. The theoretical understanding is also supported by our numerical studies.  Subsequently, we generalize these findings to settings where density functions vanish at different polynomial rates with (unknown) parameters, and establish adaptation to the H$\ddot{\text{o}}$lder smoothness parameter. 
 
 The rest of this paper is organized as follows. Section~\ref{sec:pro} introduces the problem formulation. The minimax optimal transfer learning rate is established in Section~\ref{sec:theory}. Numerical illustrations are given in Section~\ref{sec:simu}. Section~\ref{sec:ext} extends the results to more general settings. Concluding remarks are provided in Section~\ref{sec:con}. Technical proofs and lemmas are given in the Appendix and the Supplementary Materials (SM). 

\section{Problem Formulation}
\label{sec:pro}
We observe $n + n_T$ independent pairs $(X_1, Y_1), \cdots, (X_n, Y_n)$, $\cdots, (X_{n+n_T}, \cdots, Y_{n+n_T}) \in [0, 1] \times R$ with 
\begin{eqnarray*}
X_i &\sim& h, \; \; \text{for} \; i = 1, \cdots, n, \\
X_i &\sim& h_T, \; \; \text{for} \; i = n+1, \cdots, n + n_T, \\
Y_i &=& f(X_i) + \epsilon, \; \; \text{for} \; i = 1, \cdots, n+n_T, 
\end{eqnarray*}
where the error $\epsilon_i$ and $X_i$ are independent, $\epsilon_i$ independently follows a normal distribution with mean 0 and variance $\sigma^2$, and 
$\sigma \leq \overline{\sigma} < \infty$ for some $\overline{\sigma} > 0$. Let $\mcD = (X_i, Y_i)_{i=1}^n$ and $\mcD^T = (X_i, Y_i)_{i= n+1}^{n + n_T}$ denote the source and target data, respectively. Note that our theoretical results allow $n$ and $n_T$ to be 
general as long as $n + n_T \to \infty$. In particular, it is possible that $n_T \ll n$ and even $n_T = 0$.

\subsection{Minimax risks}
We measure loss of an estimate of $f$ in terms of a squared norm. Let $\| \cdot \|_{L_q(h)}$ denote the $L_q$ norm weighted by the density function $h$, i.e., for any 
$g$, 
$
\|g\|_{L_q(h)} = \left( \int |g(x)|^q h(x)dx\right)^{1/q}, q \geq 1.
$ 
Let $\mcS(\mcD)$, $\mcS(\mcD^T)$, $\mcS(\mcD, \mcD^T)$ denote the set of all estimators of $f$ based on the data set $\mcD$, $\mcD^T$, and $(\mcD, \mcD^T)$, respectively.  In this paper, we focus on the minimax rate 
under $L_2(h_T)$ loss based on source data and target data together (i.e., $(\mcD, \mcD^T)$). Given an estimator $\hat{f}_{n, n_T}$ based on $(\mcD, \mcD^T)$, let 
\[
R_{h_T}(f; \hat{f}_{n, n_T}; n, n_T) := \mbE \int_0^1 \left( f(x) - \hat{f}_{n, n_T}(x) \right)^2 h_T(x) dx
\]
be its risk, where the expectation is taken with respect to the randomness of $\mcD^T$ and $\mcD$ under the true regression function $f$. The minimax risk is 
\begin{eqnarray*}
R_{h_T}(\mcH(\beta; \kappa); \overline{\sigma}; \mcD, \mcD^T) &:=& \inf_{\hat{f}_{n, n_T} \in \mcS(\mcD, \mcD^T)}\sup_{f \in \mcH(\beta, \kappa), 0 < \sigma \leq \overline{\sigma}}R_{h_T}(f; \hat{f}_{n, n_T}; n, n_T).
\end{eqnarray*}
It describes how well one can estimate $f$ uniformly over the function class $\mcH(\beta, \kappa)$ using the data $(\mcD, \mcD^T)$ in terms of the $L_2(h_T)$ loss that is tailored towards the target population. We call $R_{h_T}(\mcH(\beta; \kappa); \overline{\sigma}; \mcD, \mcD^T)$ the transfer learning rate (TLR) based on both data sets.  For comparison, let $R_{h_T}(\mcH(\beta, \kappa); \overline{\sigma}; \mcD)$, and $R_{h_T}(\mcH(\beta, \kappa); \overline{\sigma}; \mcD^T)$ be the minimax risk of learning $f$ based on data set $\mcD$ and $\mcD^T$, respectively. According to \cite{yang1999information}, with the well-known 
$L_{\infty}$ and $L_2$ metric entropy orders of $\mcH(\beta, \kappa)$, it can be easily obtained that $R_{h_T}(\mcH(\beta, \kappa); \overline{\sigma}; \mcD^{T}) \asymp n_T^{-\frac{2\beta}{2\beta + 1}}$. However, to our knowledge, the rate of convergence of $R_{h_T}(\mcH(\beta, \kappa); \overline{\sigma}; \mcD)$ is unknown prior to this work. Note that 
$R_{h_T}(\mcH(\beta, \kappa); \overline{\sigma}; \mcD)$ is a special case of 
$R_{h_T}(\mcH(\beta; \kappa); \overline{\sigma}; \mcD, \mcD^T)$ with $n_T = 0$.

\section{Transfer Learning Rate}
\label{sec:theory}

In this Section, we explicitly derive minimax lower and upper bounds for the TLR. Throughout the rest of the paper, we consider 
$h(x) \sim Beta(a, 1)$ with $a > 0$ and $h_T(x) \sim Beta(1, 1)$, where $Beta(a_1, a_2)$ has density  $\frac{x^{a_1 - 1}(1 - x)^{a_2 - 1}}{B(a_1, a_2)} \bm{I}(x \in [0, 1])$ and $B(a_1, a_2) = \int_0^1 x^{a_1 - 1}(1 - x)^{a_2 - 1}dx$ for $a_1, a_2 > 0$. Let $a_n \succeq b_n$ denote $b_n = O(a_n)$, and $a_n \asymp b_n$ means $a_n \succeq b_n$ and $b_n \succeq a_n$. 

\subsection{Minimax lower bounding}

We give a minimax lower bound of learning $f$ based on $\mcD$ and $\mcD^T$ together. 

\begin{theorem}
\label{lem:ratef-beta-lower}
For any $a > 0$ and $\beta > 0$, the minimax learning risk satisfies:
\begin{itemize}
\item[{(I)}] when $a \leq 2 + \frac{1}{2\beta}$, 
\[ R_{h_T}(\mcH(\beta, K); \overline{\sigma}; \mcD, \mcD^T) \succeq (n + n_T)^{-\frac{2\beta}{2\beta + 1}};
\]
\item[{(II)}] when $a > 2 + \frac{1}{2\beta}$, 
\begin{eqnarray*}
R_{h_T}(\mcH(\beta, K);  \overline{\sigma}; \mcD, \mcD^T) &\succeq& 
\left\{\begin{array}{cc}
n^{-\frac{2\beta + 1}{2\beta + a}} & n_T \leq c_1n^{\frac{2\beta + 1}{2\beta + a}}\\
(\frac{n_T}{n})^{\frac{1}{a - 1}}n_T^{-\frac{2\beta}{2\beta + 1}} & c_2 n^{\frac{2\beta + 1}{2\beta + a}} \leq n_T \leq c_3n\\
n_T^{-\frac{2\beta}{2\beta + 1}} & n_T \geq c_4n.
\end{array}
\right.
\end{eqnarray*}
\end{itemize}
Here and hereafter, $c_j$ are positive constants, which are allowed to depend on $a$, $\beta$ and $\overline{\sigma}$, but are independent of  $n$ and $n_T$. Furthermore, the above bounds hold for all choices of the constants $c_1, c_2, c_3$ and $c_4$. The inequality $A_{n, n_T} \succeq B_{n, n_T}$ means that $A_{n, n_T} \geq c B_{n, n_T}$ for some positive constant $c$ that does not depend on $n$ and $n_T$, but may be dependent on $a$, $\beta$, $\overline{\sigma}$ and the constants $c_j$.

\end{theorem}

\textbf{Proof of Theorem~\ref{lem:ratef-beta-lower}}. We prove the theorem based on global entropy for Case (I) and part of Case (II), but rely on local entropy for the more challenging situations in Case (II). For minimax lower bounding, it suffices to consider a fixed $\sigma > 0$.

Case (I): We first prove that the lower rate $(n + n_T)^{-\frac{2\beta}{2\beta + 1}}$ is a valid minimax lower bound for all the situations, although it may be sub-optimal (too small) in some situations as depicted in the theorem.

Let $V_{\infty}(\epsilon)$ be a continuous upper bound of the covering $\epsilon$-entropy of a function class $\mcF$ under the sup-norm metric. Given 
$\epsilon > 0$, let $\mcG_{\epsilon} = \{f_1, f_2, \cdots, f_N\}$ be an $\epsilon$-net of $\mcF$ under the sup-norm, where 
$N \leq \exp(V_{\infty}(\epsilon))$. Let $P_f(x, y) = \frac{1}{\sqrt{2\pi}\sigma}\exp\left(-\frac{1}{2\sigma^2}( y - f(x))^2\right)h(x)$ and 
$P^T_f(x, y) = \frac{1}{\sqrt{2\pi}\sigma}\exp\left(-\frac{1}{2\sigma^2}( y - f(x))^2\right)h_T(x)$ be the joint density of 
$(X, Y)$ from the source and target distribution, respectively.

Consider any finite subset $\mcF_0 = \{\tilde{f}_1, \cdots, \tilde{f}_m\} \subset \mcF$ and a uniform prior $\Theta$ on $\mcF_0$. From 
\cite{yang1999information}, we can bound the mutual information $I(\Theta; \mcD, \mcD^T)$ in terms of the metric entropy of $\mcF$. Indeed, 
\begin{eqnarray*}
I(\Theta; \mcD, \mcD^T) &\leq& \frac{1}{m}\sum_{k=1}^m\int \prod_{i=1}^nP_{\tilde{f}_k}(x_i, y_i)\prod_{j=1}^{n_T}P^T_{
\tilde{f}_k}(x^T_i, y^T_i) \\
&&\times \log \frac{\prod_{i=1}^nP_{\tilde{f}_k}(x_i, y_i)\prod_{j=1}^{n_T}P^T_{
\tilde{f}_k}(x^T_i, y^T_i)}{Q\left( (x_i, y_i)_{i=1}^n, (x^T_j, y^T_j)_{j=1}^{n_T} \right)} dx_1dy_1 \cdots dx^T_{n_T}dy^T_{n_T},
\end{eqnarray*}
where $Q$ is taken as 
\[
Q\left( (x_i, y_i)_{i=1}^n, (x^T_j, y^T_j)_{j=1}^{n_T} \right) = \frac{1}{N}\sum_{e=1}^N\prod_{i=1}^nP_{f_e}(x_i, y_i)\prod_{j=1}^{n_T}P^T_{
f_e}(x^T_i, y^T_i).
\]
For any $\tilde{f}_k \in \mcF_0$, there exists a $f_e \in \mcG_{\epsilon}$ such that $\|\tilde{f}_k - f_e\|_{\infty} \leq \epsilon$. Then it follows that 
\begin{eqnarray*}
&&\int \prod_{i=1}^nP_{\tilde{f}_k}(x_i, y_i)\prod_{j=1}^{n_T}P^T_{
\tilde{f}_k}(x^T_i, y^T_i) \\
&& \hspace{2cm} \times \log \frac{\prod_{i=1}^nP_{\tilde{f}_k}(x_i, y_i)\prod_{j=1}^{n_T}P^T_{
\tilde{f}_k}(x^T_i, y^T_i)}{Q\left( (x_i, y_i)_{i=1}^n, (x^T_j, y^T_j)_{j=1}^{n_T} \right)} dx_1dy_1 \cdots dx^T_{n_T}dy^T_{n_T} \\
&\leq& \log N + \int \prod_{i=1}^nP_{\tilde{f}_k}(x_i, y_i)\prod_{j=1}^{n_T}P^T_{
\tilde{f}_k}(x^T_i, y^T_i) \\
&&\hspace{2cm} \times \log \frac{\prod_{i=1}^nP_{\tilde{f}_k}(x_i, y_i)\prod_{j=1}^{n_T}P^T_{
\tilde{f}_k}(x^T_i, y^T_i)}{
\prod_{i=1}^nP_{f_e}(x_i, y_i)\prod_{j=1}^{n_T}P^T_{
f_e}(x^T_i, y^T_i)
} dx_1dy_1 \cdots dx^T_{n_T}dy^T_{n_T}\\
&=&\log N + \frac{n}{2\sigma^2}\|\tilde{f}_k - f_e\|^2_h + \frac{n_T}{2\sigma^2}\|\tilde{f}_k - f_e\|^2_{h_T}\\
& \leq & \log N + \frac{1}{2\sigma^2}(n+n_T)\|\tilde{f}_k - f_e\|^2_{\infty}\\
&\leq& \log N + \frac{1}{2\sigma^2}\bar{n}\epsilon^2,
\end{eqnarray*}
where $\bar{n} = n + n_T$ is the total sample size. We then conclude that 
$I(\Theta; \mcD, \mcD^T) \leq V_{\infty}(\epsilon) + \frac{1}{2\sigma^2}\bar{n}\epsilon^2$.  
Choose $\epsilon_n$ such that $V_{\infty}(\epsilon) = \frac{1}{2\sigma^2}\bar{n}\epsilon^2$
 and let 
$\underline{\epsilon}_n$ satisfy 
\[
M_{h_T}(\underline{\epsilon}_n) = \frac{2}{\sigma^2}\bar{n}\epsilon_n^2 + 2 \log 2,
\] 
where $M_{h_T}(\epsilon)$ denotes the $\epsilon$-packing entropy of $\mcF$ under $\|\cdot\|_{h_T}$ distance. From \cite{yang1999information}, we know 
\[
\inf_{\hat{f}}\sup_{f \in \mcF, 0 < \sigma \leq \overline{\sigma}} \mbE\|\hat{f} - f\|^2_{h_T} \geq \frac{1}{8}\underline{\epsilon}_n^2,
\]
where $\hat{f}$ is based on $\mcD$ and $\mcD^T$. For $\mcF$ being the H$\ddot{\text{o}}$lder class, it is well-known \citep{carl1981entropy,temlyakov1989estimates,birman1974quantitative} that 
\[
V_{\infty}(\epsilon) \asymp M_{h_T}(\epsilon) \asymp \left(\frac{1}{\epsilon} \right)^{\frac{1}{\beta}}.
\] 
Consequently, we know $\epsilon_n \asymp \underline{\epsilon}_n \asymp \bar{n}^{-\frac{\beta}{2\beta + 1}}$. Therefore, we have shown 
\[
\inf_{\hat{f}}\sup_{f \in \mcH(\beta, \kappa), 0 < \sigma \leq \overline{\sigma}} \mbE\|\hat{f} - f\|^2_{h_T} \succeq \bar{n}^{-\frac{2\beta}{2\beta + 1}}.
\]
That is, regardless of the value of $a$, the minimax rate of convergence cannot be faster than the standard rate $\bar{n}^{-\frac{2\beta}{2\beta + 1}}$ based on 
$\bar{n}$ observations from the target distribution. 

Case (II): We derive the lower rate under different settings of $n$ and $n_T$ when $a > 2 + \frac{1}{2\beta}$.

Situation 1: $n_T \leq c_1 n^{\frac{2\beta + 1}{2\beta + a}}$ for some constant $c_1 > 0$. 

Consider the following function class $\mcF_{\eta}  = \{f \in \mcF: nD(P_f || P_{f_0}) + n_TD(P_f^T || P_{f_0}^T) \leq \eta^2\}$, for some $\eta > 0$ and $f_0 \equiv 0$.  Then for any finite set in $\mcF_{\eta}$ and a uniform prior $\Theta$ on it, by the definition of $\mcF_{\eta}$, we have that the mutual information satisfies:
\[
I(\Theta;  \mcD, \mcD^T) \leq \eta^2.
\]
Proper choices of $\eta$ will be made later. Given $\beta > 0$, let $\psi: R \to [0, \infty)$ be a $C^{\infty}(R)$ function supported in $[0, 1]$
satisfying $\psi(0) = \psi(1) = 0$, $\|\psi\|_{\infty} \leq A $ for some $ 0 < A < \infty$, $\int_0^1 \psi^2(x) dx > 0$ and $\psi \in \mcH(\beta, K)$ when restricted on $[0, 1]$. Given an integer $m > 1$, 
for $1 \leq j \leq m$, consider $\psi_j(x) = m^{-\beta} \psi(mx - (j - 1))$. Clearly, $\psi_j$ has support on $\left[\frac{j-1}{m}, \frac{j}{m}\right]$,
$\int_0^1 \psi_j^2(x)dx = \frac{\int_0^1 \psi^2(y) dy}{m^{2\beta + 1}}$, and 
\begin{eqnarray*}
a\int_0^1 \psi_j^2(x) x^{a-1}dx &=& am^{-2\beta}\int_{\frac{j-1}{m}}^{\frac{j}{m}}\psi^2(mx - (j - 1))x^{a-1}dx\\
&=& am^{-2\beta}\int_0^1\psi^2(y) \frac{ \left( y + (j-1) \right)^{a-1}}{m^{a}} dy\\
&=& \frac{a}{m^{2\beta + a}} \int_0^1 \psi^2(y) \left( y + (j-1) \right)^{a-1}dy \\
&\leq& \frac{aA^2}{m^{2\beta + a}}j^{a-1}.
\end{eqnarray*}
Now for $1 \leq J \leq m$, let $\mcT_{J} = \{-1 , 1\}^{J}$. Then for 
$t = (t(1), \cdots, t(J)) \in \mcT_{J}$, define 
\[
f_{t}(x) = \sum_{j=1}^Jt(j)\psi_j(x), \; \; \; \text{for} \; \; 0 \leq x \leq 1.
\]
Note that with a proper scaling of $\psi$, we have $f_t \in \mcH(\beta, K)$, for each $t \in \mcT_J$, 
\begin{eqnarray*}
D\left(P_{f_t} || P_{f_{0}}\right) &\leq& \frac{1}{2\sigma^2}\sum_{j=1}^{J} \frac{aA^2}{m^{2\beta + a}}j^{a-1} \leq \frac{1}{2\sigma^2}\frac{aA^2}{m^{2\beta + a}}J^{a}, \\
D\left( P^T_{f_t} || P^T_{f_{0}} \right) &=& \frac
{1}{2\sigma^2}\frac{J}{m^{2\beta + 1}} \int_0^1 \psi^2(y)dy.
\end{eqnarray*}
Therefore, we know $f_t \in \mcF_{\eta}$ for all $t \in \mcT_J$ if 
$\eta^2 \geq \frac{1}
{2\sigma^2}\left(n\frac{aA^2}{m^{2\beta + a}}J^{a} + n_T\frac{\int_0^1 \psi^2(y)dy}{m^{2\beta + 1}}  J \right)$.

We take $\eta^2 = 1$, $m = \lfloor c n^{\frac{1}{2\beta + a}} \rfloor$, and $J$ fixed. Then when $c$ is chosen large enough, we can easily verify that for all $t \in \mcT_J$, $f_{t} \in \mcF_{\eta}$. Note that the Hamming distance between $t_1$ and $t_2$ in $\mcT_J$ is at least $1$ for $t_1 \neq t_2$, and then 
\[
\int_0^1 \left( f_{t_1}(x) - f_{t_2}(x) \right)^2 dx \geq 4 \frac{\int_0^1 \psi^2(y) d y}{m^{2\beta + 1}},
\]
which is of order $n^{-\frac{2\beta + 1}{2\beta + a}}$.  When $J$ is chosen large enough ($J \geq 5$), 
the application of Fano's inequality (e.g., Equation (1), \cite{yang1999information}) readily yields the minimax lower rate $n^{-\frac{2\beta + 1}{2\beta + a}}$.

Situation 2:  
$
c_2n^{\frac{2\beta + 1}{2\beta + a}} \leq n_T \leq c_3n,
$
for some positive constants $c_2$ and $c_3$.

In this situation, we take $\eta^2 = c \left( \frac{n_T^{\frac{2\beta + a}{2\beta + 1}}}{n} \right)^{\frac{1}{a - 1}}$ for some constant $c > 0$. We choose 
$m = \lfloor c'n_T^{\frac{1}{2\beta + 1}} \rfloor$ and $J = \lfloor m \left( \frac{n_T}{n} \right)^{\frac{1}{a-1}} \rfloor$. Then with $c'$ large enough, we can ensure 
\[
\eta^2 \geq n \frac{aA^2}{m^{2\beta + a}}J^a + n_T \frac{\int_0^1 \psi^2(y)dy}{m^{2\beta + 1}} J.
\]
In $\mcT_{J}$, by Gilbert-Varshamov bound, there exist at least $2^{J/8}$ elements such that they are pairwise separated by at least 
$\frac{J}{8}$ in Hamming distance. Note that for any such pair, say $t_1$ and $t_2$, due to the construction of $\psi_j$'s, 
\[
\int_0^1 \left( f_{t_1}(x) - f_{t_2}(x) \right)^2 dx \geq 4 \frac{\int_0^1 \psi^2(y) d y}{m^{2\beta + 1}} \frac{J}{8},
\]
which is of order $n_T^{-\frac{2\beta}{2\beta + 1}} \left( \frac{n_T}{n} \right)^{\frac{1}{a-1}}$. And the log-cardinality of the set is at least $\frac{J}{8}$, which is of order 
\[
n_T^{\frac{1}{2\beta + 1}} \left( \frac{n_T}{n} \right)^{\frac{1}{a-1}} = 
\left( \frac{n_T^{\frac{2\beta + a}{2\beta + 1}}}{n} \right)^{\frac{1}{a-1}} \asymp \eta^2.
\]
With a choice of $c'$ large enough, by Fano's inequality, we conclude the minimax rate is 
\[
n_T^{-\frac{2\beta}{2\beta + 1}} \left( \frac{n_T}{n} \right)^{\frac{1}{a-1}}.
\]

Situation 3: $n_T \geq c_4n$ for some constant $c_4 > 0$.  The lower rate is readily seen in the proof of Case (I).

Based on all above, Theorem~\ref{lem:ratef-beta-lower} follows.

\subsection{Minimax upper bounding}
To obtain a good upper bound, adopting a technique of \cite{schmidt2022local}, we consider a local region measured by a function $t_n(x)$ that satisfies 
\begin{eqnarray}
\label{eq:lb}
t^{2\beta}_n(x) = \frac{1}{nH(x \pm t_n(x)) + n_TH_T(x \pm t_n(x))},
\end{eqnarray}
where $H(x \pm t_n(x)) = \int_{\max(x - t_n(x), 0)}^{\min(x + t_n(x), 1)}h(y)dy =(\min(x + t_n(x), 1))^a - (\max(x - t_n(x), 0))^a$ and $H_T(x \pm t_n(x)) = \int_{\max(x - t_n(x), 0)}^{\min(x + t_n(x), 1)}h_T(y)dy = \min(x + t_n(x), 1) - \max(x - t_n(x), 0)$. 
Clearly, when local polynomial tools \citep{fan2018local} are used to estimate the nonparametric function $f(x)$, the function $t_n(x)$ can be regarded as a local (or adaptive) bandwidth depending on the point $x$. Moreover, based on local polynomial techniques, it is known that the bias and variance of nonparametric estimator are at the order $t_n^{\beta}(x)$ and $(nt_n(x))^{-1}$, respectively. Thus, equation~\eqref{eq:lb} is to provide a good trade off between the bias and variance.  Following similar arguments as those in \cite{schmidt2022local}, Lemma~\ref{lem:tn} and Lemma~\ref{lem:tn1} in the Appendix show the existence of $t_n(x)$ and several properties of $t_n(x)$ including the uniqueness of $t_n(x)$ and its differentiality. We also refer to $t_n(x)$ as the spread function when it satisfies equation~\eqref{eq:lb}. The next Theorem~\ref{lem:ratef-beta} gives an upper bound for the function class $\mcH(\beta, \kappa)$ based on local polynomial estimators, defined in equation~\eqref{eq:lp}, when $\beta > 1$ and Nadaraya-Watson (NW) estimators, defined in equation~\eqref{eq:nw}, when $\beta \in (0, 1]$.

\begin{theorem}
\label{lem:ratef-beta}
For any $a > 0$ and $\beta > 0$, the minimax learning rate satisfies:  
\begin{itemize}
\item[{(I)}] when $a \leq 2 + \frac{1}{2\beta}$, 
\[
R_{h_T}(\mcH(\beta, K); \overline{\sigma}; \mcD, \mcD^T) \preceq (n + n_T)^{-\frac{2\beta}{2\beta + 1}};
\]
\item[{(II)}] when $a > 2 + \frac{1}{2\beta}$, 
\begin{eqnarray*}
R_{h_T}(\mcH(\beta, K); \overline{\sigma}; \mcD, \mcD^T) &\preceq& 
\left\{
\begin{array}{cc}
n^{-\frac{2\beta + 1}{2\beta + a}}  & n_T \leq c_1n^{\frac{2\beta + 1}{2\beta + a}}\\
(\frac{n_T}{n})^{\frac{1}{a - 1}}n_T^{-\frac{2\beta}{2\beta + 1}} & c_2n^{\frac{2\beta + 1}{2\beta + a}} \leq n_T \leq c_3n\\
n_T^{-\frac{2\beta}{2\beta + 1}} & n_T \geq c_4n.
\end{array}
\right.
\end{eqnarray*}
\end{itemize}
\end{theorem}

\textbf{Proof of Theorem~\ref{lem:ratef-beta}}. Using   Lemma~\ref{lem:upb} (NW estimator for $\beta \in (0, 1]$) and Lemma~\ref{local-p} (local polynomial estimator for $\beta > 1$), we have that 
\[
\mbE\left[ (\hat{f}_{n, n_T}(x) - f(x))^2 \right] \leq C \left( (t_n(x))^{2\beta} + \frac{1}{nH(x \pm t_n(x)) + n_TH_T(x \pm t_n(x))} \right),
\]
where $C$ is a constant independent of $x$ for $x \in [0, 1]$. Based on the definition of $t_n(x)$, we have that  
\begin{eqnarray}
\label{eq:mix}
\int_{0}^1\left( \hat{f}_{n, n_T}(x) - f(x)\right)^2 h_T(x) dx \preceq \int_0^1 (t_n(x))^{2\beta} dx .
\end{eqnarray}
Next we show the order of the spread function $t_n(x)$. We split $[0, 1]$ into three intervals 
$I_1 := [0, \alpha_n]$, $I_2 := (\alpha_n,  b_n)$ and $I_3 := [b_n, 1]$ such that $t_n(\alpha_n) = \alpha_n$ and $t_n(b_n) = 1 - b_n$.

The solution $\alpha_n$ of the equation $t_n(x) = x$ must satisfy 
\[
\alpha_n^{2\beta} = \frac{1}{n (2\alpha_n)^{a} + n_T(2\alpha_n)}.
\]
Let $a_{n, n_T} = o(b_{n, n_T})$ represents $a_{n, n_T}/b_{n, n_T} \to 0$ as $n, n_T \to \infty$. Clearly, if
$n_T2\alpha_n^{2\beta + 1} = o\left(n2^a\alpha_n^{2\beta + a}\right)$ and $n_T = o\left(n^{\frac{2\beta + 1}{2\beta + a}}\right)$, it follows that $\alpha_n \asymp n^{-\frac{1}{2\beta + a}}$; when 
$n2^{a - 1}\alpha_n^{2\beta + a} = o\left(n_T\alpha_n^{2\beta + 1}\right)$ and $n^{\frac{2\beta + 1}{2\beta + a}} = o\left(n_T\right)$, it follows that $\alpha_n \asymp n_T^{-\frac{1}{2\beta + 1}}$. Combining the results, we have $\alpha_n \asymp \left( n^{\frac{2\beta + 1}{2\beta + a}} + n_T \right)^{-\frac{1}{2\beta + 1}}$. 
Note that $t_n(0)^{2\beta} = \frac{1}{n(t_n(0))^{\alpha} + n_Tt_n(0)}$. Using similar arguments, we have that $t_n(0) \asymp \left( n^{\frac{2\beta + 1}{2\beta + a}} + n_T \right)^{-\frac{1}{2\beta + 1}}$.  According to Lemma 2 (iii), we know that $t_n(x)$ is strictly decreasing in the region $[0, b_n)$, and strictly increasing in the region $[b_n, 1]$ and thus $t_n(\alpha_n) \leq t_n(x) \leq t_n(0)$ for $x \in [0, \alpha_n]$, that is, 
\[
t_n(x) \asymp \left( n^{\frac{2\beta + 1}{2\beta + a}} + n_T \right)^{-\frac{1}{2\beta + 1}}, \; \; x \in \left[0,  \alpha_n\right].
\] 
Similarly, for $b_n$ we can obtain that $1 - b_n \asymp \left(n + n_T \right)^{-\frac{1}{2\beta + 1}}$, 
$t_n(1) \asymp \left(n + n_T \right)^{-\frac{1}{2\beta + 1}}$, and $t_n(b_n) \leq t_n(x) \leq t_n(1)$ for $x \in [b_n, 1]$. That is, 
\[
t_n(x) \asymp \left( n + n_T \right)^{-\frac{1}{2\beta + 1}}, \; \; x \in \left[b_n,  1\right].
\]
For the second region, we can obtain that 
$t_n(x) \asymp \left( \frac{1}{nx^{a - 1} + n_T } \right)^{\frac{1}{2\beta + 1}}$.

Then, we have 
\begin{eqnarray*}
\mbE \int_0^1 \big( \hat{f}_{n, n_T}(x) - f(x) \big)^2 h_T(x) dx 
&\preceq& \int_0^1 \big( t_n(x) \big)^{2\beta}dx \\
&\preceq&  \left( n^{\frac{2\beta + 1}{2\beta + a}} + n_T \right)^{-1}  + \int_{\alpha_n}^{1}  \left( \frac{1}{nx^{a - 1} + n_T } \right)^{\frac{2\beta}{2\beta + 1}}dx.
\end{eqnarray*}

Next, we give the explicit form of the upper bound under different settings of $n$ and $n_T$ for Cases (I)--(II):

Situation 1: If $n_T \leq c_1n^{\frac{2\beta + 1}{2\beta + a}}$, then it follows that 
\begin{eqnarray*}
\mbE \int_0^1 \big( \hat{f}_{n, n_T}(x) - f(x) \big)^2 h_T(x) dx &\preceq& n^{-\frac{2\beta + 1}{2\beta+ a}} + 
n^{-\frac{2\beta}{2\beta + 1}}\int_{\alpha_n}^{1} x^{(1 - a)\frac{2\beta}{2\beta + 1}}dx \\
&\preceq& n^{-\frac{2\beta + 1}{2\beta+ a}} + n^{-\frac{2\beta}{2\beta + 1}}(\log n)^{\bm{I}(a = 2 + \frac{1}{2\beta})} \\
&\preceq& \left\{ \begin{array}{cc}
n^{-\frac{2\beta}{2\beta + 1}}, & \text{for} \; 0 < a < 2 + \frac{1}{2\beta}\\
n^{-\frac{2\beta}{2\beta + 1}} \log n, & \text{for} \; a = 2 + \frac{1}{2\beta}\\
n^{-\frac{2\beta + 1}{2\beta+ a}}, & \text{for} \; a > 2 + \frac{1}{2\beta}.\\
\end{array}\right.
\end{eqnarray*}

Situation 2: If $a > 1$ and $c_2n^{\frac{2\beta + 1}{2\beta + a}} \leq n_T \leq c_3n$, then it follows from $
b_n \asymp 1 - n^{-\frac{1}{2\beta+1}} < 1$ that 
\begin{eqnarray*}
&&\mbE \int_0^1 \big( \hat{f}_{n, n_T}(x) - f(x) \big)^2 h_T(x) dx \\
&\preceq& n_T^{-1} + \int_{\left(\frac{1}{n_T}\right)^{\frac{1}{2\beta + 1}}}^{1}\left( \frac{1}{nx^{a - 1} + n_T} \right)^{\frac{2\beta}{2\beta + 1}}dx\\
&\preceq&n_T^{-1} + \int_{\left(\frac{1}{n_T}\right)^{\frac{1}{2\beta + 1}}}^{(\frac{n_T}{n})^{1/(a - 1)}} n_T^{-\frac{2\beta}{2\beta + 1}}dx + n^{-\frac{2\beta}{2\beta + 1}}\int_{(\frac{n_T}{n})^{1/(a - 1)}}^{1}
x^{(1 - a)\frac{2\beta}{2\beta + 1}}dx\\
&\preceq&\left\{
\begin{array}{cc}
n_T^{-1} + (\frac{n_T}{n})^{1/(a - 1)} n_T^{-\frac{2\beta}{2\beta + 1}} + n^{-\frac{2\beta}{2\beta + 1}}, & \text{for} \; 1 < a < 2 + \frac{1}{2\beta}\\
n_T^{-1} + (\frac{n_T}{n})^{1/(a - 1)} n_T^{-\frac{2\beta}{2\beta + 1}} + n^{-\frac{2\beta}{2\beta + 1}}\log n, & \text{for} \;  a = 2 + \frac{1}{2\beta}\\
n_T^{-1} + 2(\frac{n_T}{n})^{1/(a - 1)} n_T^{-\frac{2\beta}{2\beta + 1}} - n^{-\frac{2\beta}{2\beta + 1}}, 
& \text{for} \; a > 2 + \frac{1}{2\beta}\\
\end{array}
\right.\\
&\preceq&\left\{
\begin{array}{cc}
n^{-\frac{2\beta}{2\beta + 1}}, & \text{for} \; 1 < a < 2 + \frac{1}{2\beta}\\
  n^{-\frac{2\beta}{2\beta + 1}}\log n, & \text{for} \; a = 2 + \frac{1}{2\beta}\\
(\frac{n_T}{n})^{1/(a - 1)} n_T^{-\frac{2\beta}{2\beta + 1}}, & \text{for} \; a > 2 + \frac{1}{2\beta}.\\
\end{array}
\right.
\end{eqnarray*}

Situation 3: If $n_T \geq c_{4}\left( n I (a > 1) + n^{\frac{2\beta + 1}{2\beta + a}}I(a \leq 1)\right)$, it follows that 
\begin{eqnarray*}
\mbE \int_0^1 \big( \hat{f}_{n, n_T}(x) - f(x) \big)^2 h_T(x) dx &\preceq& n_T^{-1} + 
\int_{\left(\frac{1}{n_T}\right)^{\frac{1}{2\beta + 1}}}^{1} n_T^{-\frac{2\beta}{2\beta + 1}}dx \\
&\preceq& n_T^{-1}+ n_T^{-\frac{2\beta}{2\beta + 1}} \preceq n_T^{-\frac{2\beta}{2\beta + 1}}.
\end{eqnarray*}

Noting that in Situations 1 and 2, when $a = 2 + \frac{1}{2\beta}$, $n_T \leq n$, it is sufficient to use only the source data to construct the estimator, which yields 
$\mbE \int_0^1 \big( \hat{f}_{n, n_T}(x) - f(x) \big)^2 h_T(x) dx \preceq n^{-\frac{2\beta}{2\beta + 1}}$. Then, Theorem~\ref{lem:ratef-beta} follows. 

\textbf{Remark}: \cite{schmidt2022local} obtained a similar upper bound (in proabability with an extra logarithm term needed for their sup-norm loss) for a special case with $\beta = 1$ and $a > 2.5$. Particularly, they showed that (Lemma 17), under the condition $n^{\frac{3}{2+a}}\log^{\frac{a-1}{2 + a}}n \ll n_T < n$, 
 \[
 \int_0^1 (\hat{f}_{n, n_T}(x) - f_0(x))^2 dx \preceq \left( \frac{\log n}{n_T} \right)^{2/3}\left(\frac{n_T}{n}\right)^{1/(a - 1)},
 \]
 with probability converging to one as $n \to \infty$.

 \subsection{Optimal transfer learning rate and synergistic learning phenomenon}
Combining Theorems~\ref{lem:ratef-beta-lower} and~\ref{lem:ratef-beta}, we can obtain the optimal TLR, which is summarized in the following Corollary~\ref{th:TLR}. To show the SLP, we define the synergistic acceleration rate (SAR) as the rate of 
 \[
 S(\mcD, \mcD^T; n, n_T) = \frac{\min\left( R_{h_T}(\mcH(\beta, K); \overline{\sigma}; \mcD^T), 
R_{h_T}(\mcH(\beta, K); \overline{\sigma}; \mcD) \right)}{R_{h_T}(\mcH(\beta, K); \overline{\sigma}; \mcD, \mcD^T)},
 \]
 which is the ratio between the better rate of learning $f$ based on either the target data or the source data and the optimal TLR based on the target and source data together. 
 \begin{corollary}
 \label{th:TLR}
 For any $a > 0$ and $\beta > 0$, 
 \begin{itemize}
 \item[{(a)}] the optimal TLR is:
 \begin{eqnarray*}
R_{h_T}(\mcH(\beta, K); \overline{\sigma}; \mcD, \mcD^T) &\asymp& 
\left\{
\begin{array}{ccc}
(n + n_T)^{-\frac{2\beta}{2\beta + 1}} & \forall \; n, n_T & \; a \leq 2 + \frac{1}{2\beta} \\
n^{-\frac{2\beta + 1}{2\beta + a}}  & n_T \leq c_0n^{\frac{2\beta + 1}{2\beta + a}} & \; a > 2 + \frac{1}{2\beta}\\
(\frac{n_T}{n})^{\frac{1}{a - 1}}n_T^{-\frac{2\beta}{2\beta + 1}} & c_1n^{\frac{2\beta + 1}{2\beta + a}} \leq n_T \leq c_2n & \;  a > 2 + \frac{1}{2\beta}\\
n_T^{-\frac{2\beta}{2\beta + 1}} & n_T \geq c_3n & \; a > 2 + \frac{1}{2\beta}.
\end{array}
\right.
\end{eqnarray*}
 \item[{(b)}] For the SAR, we have 
 \begin{eqnarray*}
  S(\mcD, \mcD^T; n, n_T)  \asymp
  \left\{ 
\begin{array}{ccc}
\left( n_T n^{-\frac{2\beta + 1}{2\beta + a}} \right)^{\frac{2\beta}{2\beta + 1} - \frac{1}{a-1}}, & \; \text{for} \; \; n^{\frac{2\beta + 1}{2\beta + a}} \prec n_T \preceq n^{\frac{(2\beta + 1)^2}{2\beta (2\beta + a)}}, a > 2 + \frac{1}{2\beta}\\
\left(n/n_T\right)^{\frac{1}{a-1}} & \; \text{for} \; \; n^{\frac{(2\beta + 1)^2}{2\beta (2\beta + a)}} \preceq n_T \prec n, a > 2 + \frac{1}{2\beta}\\
1 & \text{otherwise}.\\
\end{array}
\right.  
 \end{eqnarray*}
 \end{itemize}
 \end{corollary}
 
 According to Corollary~\ref{th:TLR}, 
 we see that only when  $n^{\frac{2\beta + 1}{2\beta + a}} \prec n_T \prec n$ and $a > 2 + \frac{1}{2\beta}$, the SAR tends to infinity as $n \to \infty$, i.e., the SLP occurs. Furthermore, we can see that 
 as $a$ increases, a smaller $n_T$ is enough to be helpful on learning $f$ as long as $n_T \succ n^{\frac{2\beta + 1}{2\beta + a}} $. On the other hand, as $\beta$ increases, a larger $n_T$ is required to make the additional target data helpful for learning $f$. In the next subsection, we discuss in more details on how the SAR changes according to $n$ and $n_T$ in relation to $\beta$ and $a$.  
 
  \subsection{Sample size relationship regions and SAR}

When the SLP occurs, there are three critical target sample sizes relative to the source sample size 
that determine the behavior of the SAR: the lower critical target size 
$n_T \asymp n^{\frac{2\beta + 1}{2\beta + a}}$, after which the target sample becomes 
useful, the equalizing critical target size $n_T \asymp n^{\frac{(2\beta + 1)^2}{2\beta (2\beta + a)}}$, 
at which the source sample and the target sample give exactly the same rate of convergence on their own, and the upper 
target size $n_T \asymp n$, after which the source sample does not benefit the rate of convergence when combined with the target 
data.
 
It is interesting to note that the SAR behaves differently before and after the 
equalizing critical target sample size, i.e., $n_T \asymp n^{\frac{(2\beta + 1)^2}{2\beta (2\beta + a)}}$, where SAR achieves the largest value.
Particularly, the SAR increases as $n_T$ increases from the lower critical target size to the equalizing critical target size, and then decreases as $n_T$ keeps increasing from the equalizing critical target size to the upper critical target size, i.e., $n_T \asymp n$. Accordingly, when $a > 2 + \frac{1}{2\beta}$, we divide the sample size relationship between $n_T$ and $n$ (SSR) into four regions where source and target data serve rather different roles in the contribution to the optimal TLR. 

\begin{enumerate}
\item Source Dominating region (SD): the target size is no larger than the lower critical target size $n_T \preceq n^{\frac{2\beta + 1}{2\beta + a}}$. The need of estimation of $f$ around $x_0 = 0$, where the likelihood ratio $\rho(x) \to \infty$ as $x \to x_0$, indicates the issue that there are relatively few source observations to help learn $f$ around $x_0$. Due to the gap between $h$ and $h_T$ around $x_0$, the learning rate based on the source data only is slow at $n^{-\frac{2\beta + 1}{2\beta + a}}$, slower than the rate $n^{-\frac{2\beta}{2\beta + 1}}$ under $L_2(h)$ loss. In region SD, $n_T$ is too small to be helpful on learning $f$ around $x_0$. Thus, the learning rate is determined by the source data only.

\item Source Leading region (SL): the target size is larger than the lower critical target size and smaller than the equalizing critical target size, i.e., $n^{\frac{2\beta + 1}{2\beta + a}} \prec n_T \preceq n^{\frac{(2\beta + 1)^2}{2\beta (2\beta + a)}}$. In region SL, the learning rate based on the source data only is still faster than that based on the target data only, i.e., $n_T^{-\frac{2\beta}{2\beta + 1}}$. Combining with the $n_T$ target observations, the learning rate around $x_0$ gets improved over that based on the source data only. 
From the formula of SAR, it is easy to see that the SAR increases as $a$ increases in SL. This is intuitively clear since when the source data is leading, the supporting effect of the target data is more significant when $a$ is large. 

\item Target Leading region (TL): the target size is no smaller than the equalizing critical target size and smaller than the upper critical target size, i.e., $n^{\frac{(2\beta + 1)^2}{2\beta (2\beta + a)}} \preceq n_T \prec n$. In this region, the learning rate based on the target data only is already faster than that based on the source data only. The source data help the estimation at majority of the points where $h(x)$ is not very small due to its large size than $n_T$. Combining with $n$ source observations, the learning rate gets improved over that based on the target data only. In TL, we can see that the SAR actually decreases as $a$ increases. Clearly, when target data is leading, the supporting effect of the source data is less significant when $a$ is larger.  

\item Target Dominating region (TD): the target size is no smaller than the upper critical target size, i.e., $n_T \succeq n$. In this region, $n_T$ is  large enough such that the $n$ source observations are not useful in terms of rate of convergence.

\end{enumerate}

Table~\ref{tab:ts} displays the optimal TLR and the better individual learning rate based on the source or target data alone in all possible situations 
in our context. 

\begin{table}
\centering
\caption{The minimax optimal TLR based on both $n_T$ target observations and $n$ source observations and the better individual learning rate based on the target data only or the source data only for the H$\ddot{\text{o}}$lder class $\mcH(\beta, \kappa)$ with 
$h_T\sim \text{Beta}(1, 1), h\sim \text{Beta}(a, 1)$ for $a > 0$.}
\label{tab:ts}
\resizebox{0.94\textwidth}{!}{
\begin{tabular}{lccccc}
\hline\hline
Settings&&Region&\multicolumn{3}{c}{Learning rate} \\
&&&Optimal TLR &Better individual rate&Synergistic \\
\hline
\multirow{2}*{$0 < a \leq 2 + \frac{1}{2\beta}$}&$n_T \preceq n$&SD& $n^{-\frac{2\beta}{2\beta + 1}}$ &  $n^{-\frac{2\beta}{2\beta + 1}}$ & No\\
&$n_T \succeq n$&TD& $n_T^{-\frac{2\beta}{2\beta + 1}}$ &  $n_T^{-\frac{2\beta}{2\beta + 1}}$ & No\\
\\
\multirow{4}*{$a > 2 + \frac{1}{2\beta}$}&$n_T \preceq n^{\frac{2\beta + 1}{2\beta + a}} $ &SD&$n^{-\frac{2\beta + 1}{2\beta + a}}$&$n^{-\frac{2\beta + 1}{2\beta + a}}$ &No\\
&$n^{\frac{2\beta + 1}{2\beta + a}} \prec n_T \preceq n^{\frac{(2\beta + 1)^2}{2\beta (2\beta + a)}}$ &SL&
$\left(\frac{n_T}{n}\right)^{\frac{1}{a-1}}n_T^{-\frac{2\beta}{2\beta + 1}}$&$n^{-\frac{2\beta + 1}{2\beta + a}}$& Yes\\
&$n^{\frac{(2\beta + 1)^2}{2\beta (2\beta + a)}} \preceq n_T \prec n$&TL&$\left(\frac{n_T}{n}\right)^{\frac{1}{a-1}}n_T^{-\frac{2\beta}{2\beta + 1}}$&$n_T^{-\frac{2\beta}{2\beta + 1}}$ &Yes\\
& $n_T \succeq n$&TD&$n_T^{-\frac{2\beta}{2\beta + 1}}$&$n_T^{-\frac{2\beta}{2\beta + 1}}$  &No\\
\hline \hline
\end{tabular}
}

\end{table}

\section{Simulation Studies}
\label{sec:simu}

In this section, we conduct simulations to gain a numerical understanding of 
the main theoretical results. In particular, under the setting   
$a > 2 + \frac{1}{2\beta}$, we show numerical results in the four SSR regions to illustrate the SLP with $\beta = 0.5$ and $0.9$, respectively. Considering $a = 4$, according to Corollary~\ref{th:TLR}, we know that 
the SLP happens when $n^{\frac{2\beta + 1}{2\beta + 4}} \prec n_T \prec n$. Thus, considering $n = 3000, 5000, 10000, 30000, 50000, 100000$, we set $n_T = 0.1 \times n^{\frac{2\beta + 1}{2 \beta + 4}}, 0.5 \times n^{\frac{2\beta + 1}{2 \beta + 4}(1+ c_{4, SL})}, n^{\frac{2\beta + 1}{2 \beta + 4}(1+ c_{4, TL})}, 10\times n$, respectively, corresponding to regions SD, SL, TL, TD, respectively, with $c_{4, SL}$ and $c_{4, TL}$ properly chosen.

According to Corollary~\ref{th:TLR}, in our setting, when $a > 2 + \frac{1}{2\beta}$, we have in the region SL: 
\begin{eqnarray}
\label{eq:SL}
\log \text{SAR} =  O(1) + c_{a, SL} \times \frac{2\beta + 1}{2\beta + a} \times  \left( \frac{2\beta}{2\beta + 1} - \frac{1}{a-1} \right)\log n, 
\end{eqnarray}
 and in the region TL: 
\begin{eqnarray}
\label{eq:TL}
\log \text{SAR} = O(1) + \left(1 - \frac{2\beta + 1}{2\beta + a} \times (1 + c_{a, TL})\right) \times \frac{1}{a-1} \log n,
\end{eqnarray}
and for the SD and TD regions, $\log \text{SAR} = O(1)$.  Note that for $a \leq 2 + \frac{1}{2\beta}$, we always have $\log \text{SAR} = O(1)$.  Thus, our simulations focus on examining that (a) the linear relationship between $\log \text{SAR}$ and $\log n$ in the regions SL and TL when $a > 2 + \frac{1}{2\beta}$, (b) bounded relationship in regions SD and TD when $a > 2 + \frac{1}{2\beta}$, and (c) bounded relationship when $a < 2 + \frac{1}{2\beta}$.   

For $\beta = 0.5$, we set $(c_{4, SL}, c_{4, TL}) = (0.95, 1.2)$ representing the SL and TL regions, where the SLP should occur in the two different ways as explained earlier. Similarly, for $\beta = 0.9$, we set $(c_{4, SL}, c_{4, TL}) = (0.4, 0.7)$. For comparison, we also add the results for $(a, \beta) = (2.6, 0.5)$ and $(2.1, 0.9)$ under the same SSRs, i.e., SSR-1 with 
$n_T \preceq n^{\frac{2\beta + 1}{2\beta + 4}}$, SSR-2 with 
$n^{\frac{2\beta + 1}{2\beta + 4}} \prec n_T \preceq n^{\frac{(2\beta + 1)^2}{2\beta(2\beta + 4)}}$, SSR-3 with $n^{\frac{(2\beta + 1)^2}{2\beta(2\beta + 4)}} \preceq n_T \prec n$ and SSR-4 with $n_T \succeq n$, respectively. Clearly, for the two settings $(a, \beta) = (2.6, 0.5)$ and $(2.1, 0.9)$, the SLP does not happen. 

Two cases of $f(x)$ are considered: $f_1(x) = 0.2 \times x^{\beta}$ and $f_2(x) = 0.2 \times |x - 0.5|^{\beta}$. The data are generated from $y = f(x) + \epsilon$, where $\epsilon \sim N(0, 0.3^2)$ and $x$ either follows a uniform distribution, i.e., $h_T(x) \sim U(0, 1)$ for target data, or $x$ follows a beta distribution, i.e., $h(x) \sim Beta(a, 1)$ as defined in Section~\ref{sec:theory} for source data. 
The estimation of $f$ is calculated based on NW estimator, defined in equation~\eqref{eq:nw}, with the indicator function as the kernel function, i.e., 
$I(x - h \leq x_i \leq x + h)$, and the bandwidth $h = ct_n(x) = c \left[ \alpha_n I( 0 \leq x \leq \alpha_n) +  \left(n x^{a - 1} + n_T\right)^{-\frac{1}{2\beta + 1}}I( \alpha_n < x < b_n) + 
(1 - b_n) \right.$ 
$\left. \times I (b_n \leq x \leq  1)\right]$, where  $\alpha_n  = (n^{\frac{2\beta + 1}{2\beta + a}} + n_T)^{-\frac{1}{2\beta + 1}} , b_n = 1 - (n + n_T)^{-\frac{1}{2\beta + 1}}$ and $c$ is a constant and set to be $0.7$ in this simulation.
Denote $\hat{f}_{n, n_T}(x)$ as the NW estimator based on source data and target data together. Let $\hat{f}(x)$ and $\hat{f}_T(x)$ represent the NW estimator based on the source data (treating $n_T$ = 0) and the target data (treating $n$ = 0), respectively. 

The squared $L_2(h_T)$ loss is approximated at $N$ fixed grid points in $[0, 1]$. That is, for an estimator $\hat{f}$, we calculate  $L_{h_T}(f, \hat{f}^{(k)}) = \frac{1}{N}\sum_{i=1}^{N} (\hat{f}^{(k)}(x_i) - f(x_i))^2$, where $N = 3000$, $x_i = i/N$, and $k$ represents the $k$th simulation run.  Let 
$L^{(k)}, L_T^{(k)}$ and $L_{n, n_T}^{(k)}$ denote the loss of $\hat{f}$, $\hat{f}_T$ and $\hat{f}_{n, n_T}$ respectively based on the $k$th simulation run. For each simulation setting with the chosen $n$, $n_T$, $\beta$ and $a$, to obtain a stable risk ratio, we further calculate the median ratio in terms of $\log \text{SAR} = \log \frac{
\min \left(\text{median}_{k=1,\cdots, 100}L^{(k)}, \text{median}_{k=1, \cdots, 100}L^{(k)}_T\right)}{\text{median}_{k=1, \cdots, 100}L^{(k)}_{n, n_T}}$ to see how the $\log \text{SAR}$ changes in $\log n$ according to varying $n$, $n_T$, $\beta$ and $a$.

Based on equations~\eqref{eq:SL} and~\eqref{eq:TL},  we have a linear relationship
with slope $\{0.063, 0.04\}$ at $\beta = 0.5$ for regions SL and TL, respectively, with slope $0.060$ at $\beta = 0.9$ for both regions SL and TL.
We show the $\log \text{SAR}$ v.s. $\log n$ in the following Figures~\ref{fig-x0beta05-log} and~\ref{fig-x0beta09-log}. 
It can be seen that, when $\beta \in \{0.5, 0.9\}$, $a = 4$,  clear line shapes (linear regression with significant coefficients close to the true theoretical values) are shown for regions SL and TL, while in regions SD and TD, bounded relationships are shown, which support the results of Corollary~\ref{th:TLR} that only when $n^{\frac{2\beta + 1}{2\beta + 4}} \prec n_T \prec n$, the SLP occurs. 
In addition, with $(\beta, a) = (0.5, 2.6)$ and $(0.9, 2.1)$, bounded relationships are shown in Figures~\ref{fig-x0beta05-log} and~\ref{fig-x0beta09-log} in all four regions SSR-r, r = 1, 2, 3, 4, respectively, which also support that only when $a > 2 + \frac{1}{2\beta}$, the SLP may happen.

\begin{figure}
\centering
\includegraphics[height=14cm, width=14.5cm]{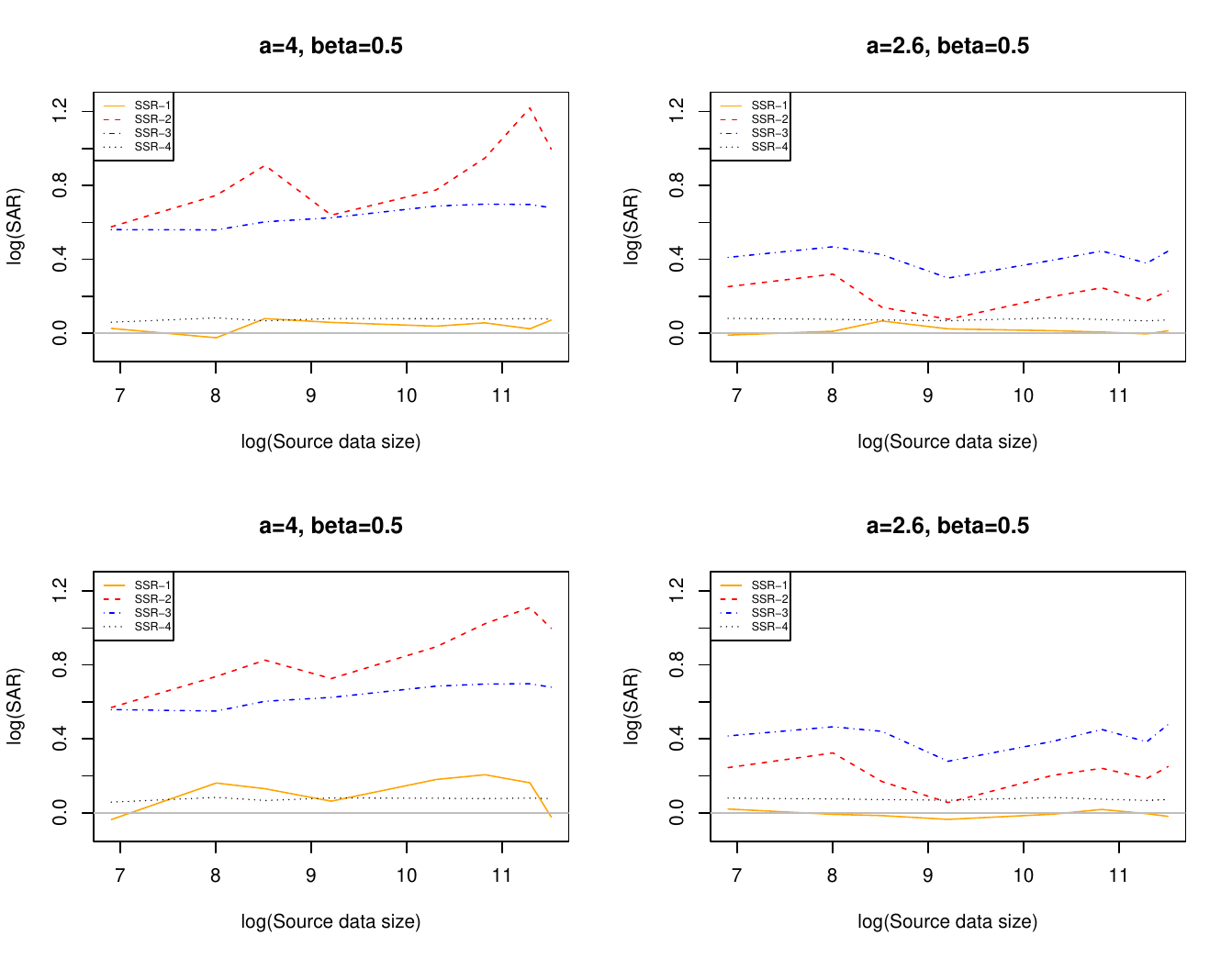}
\caption{The  logarithms of approximated SAR are plotted as $\log n$ increases, where the solid orange line represents the ratio calculated with $n_T = 0.1 \times n^{\frac{2\times 0.5 + 1}{2 \times 0.5 + 4}}$ (SSR-1), 
the dotted black line represents the ratio calculated with $n_T = 10\times n$ (SSR-4), 
the dashed red line represents the ratio calculated with $n_T = 0.5 \times n^{\frac{2\times 0.5 + 1}{2 \times 0.5 + 4}(1 + 0.95)}$ (SSR-2), and 
the dot-dashed blue line represents the ratio calculated with $n_T = n^{\frac{2\times 0.5 + 1}{2 \times 0.5 + 4}(1 + 1.2)}$ (SSR-3). In the first row, the true 
$f(x) = f_1(x) := 0.2x^{0.5}$. In the second row, the true $f(x) = f_2(x) := 0.2 |x - 0.5|^{0.5}$.}
\label{fig-x0beta05-log}
\end{figure}

\begin{figure}
\centering
\includegraphics[height=14cm, width=14.5cm]{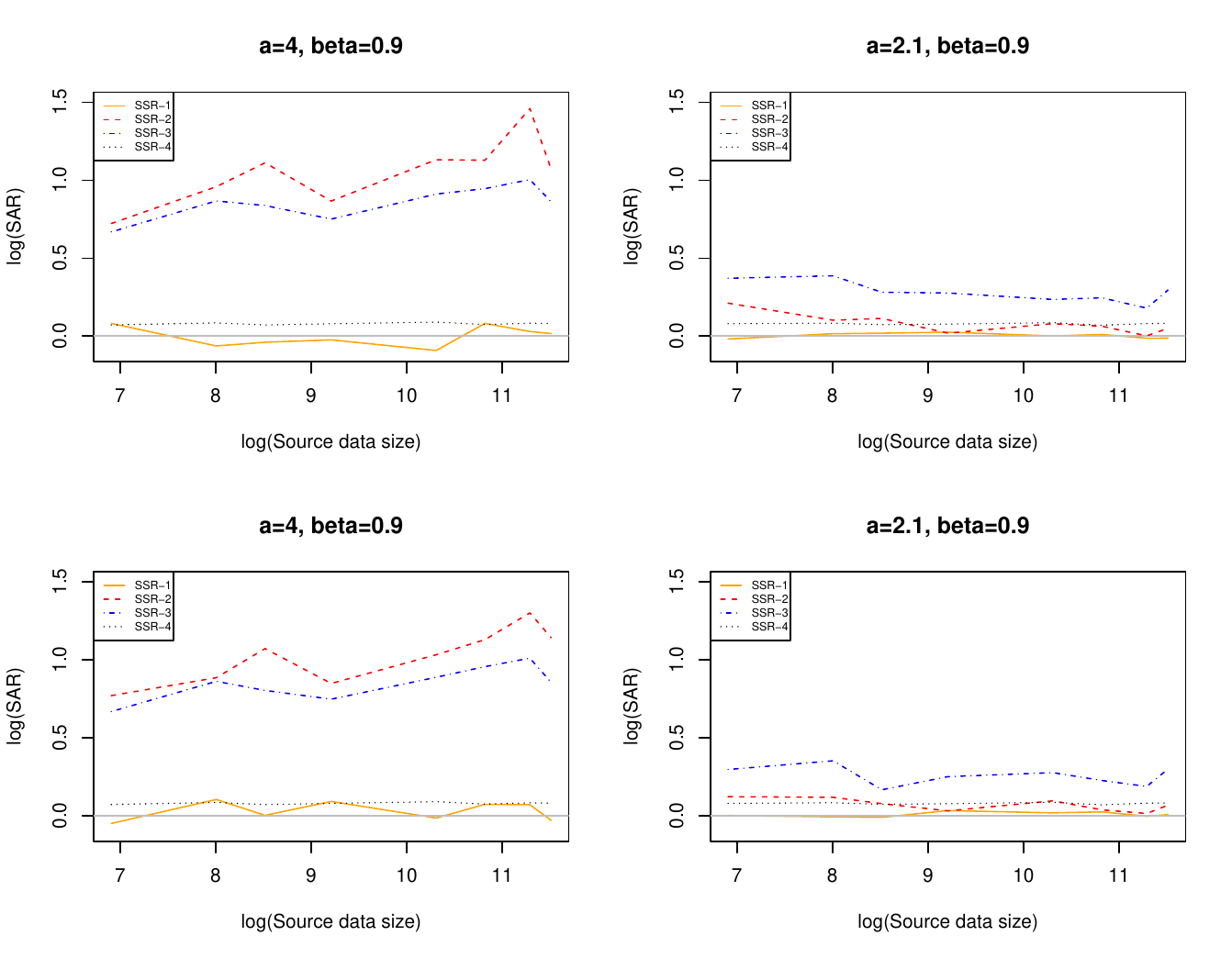}
\caption{The  logarithms of approximated SAR are plotted as $\log n$ increases, where the solid orange line represents the ratio calculated with $n_T = 0.1 \times n^{\frac{2\times 0.9 + 1}{2 \times 0.9 + 4}}$ (SSR-1), 
the dotted black line represents the ratio calculated with $n_T = 10\times n$ (SSR-4), 
the dashed red line represents the ratio calculated with $n_T = 0.5 \times n^{\frac{2\times 0.9 + 1}{2 \times 0.9 + 4}(1 + 0.4)}$ (SSR-2), and 
the dot-dashed blue line represents the ratio calculated with $n_T = n^{\frac{2\times 0.9 + 1}{2 \times 0.9 + 4}(1 + 0.7)}$ (SSR-3). In the first row, the true 
$f(x) = f_1(x) := 0.2x^{0.9}$. In the second row, the true $f(x) = f_2(x) := 0.2 |x - 0.5|^{0.9}$.}
\label{fig-x0beta09-log}
\end{figure}

\section{Extension}
\label{sec:ext}
In this section, we extend the results to more general situations. Importantly, we allow the target density and the source 
density to be general and interestingly, it turns out that the minimax transfer learning rate depends on the smoothness of 
$f$ and $h$ and $h_T$ in a more complicated way.  For notational ease in some later expressions, $h$ will be denoted as $h_A$.

\subsection{Density functions that may approach 0 at different orders}

Now, we consider a general case under which both $h_T$ and $h_A$ are upper-bounded but may approach $0$ at a finite number of points in $[0, 1]$.  Define the ``potential singularity points" to be 
$$S^o = \{x \in [0, 1]: {\lim}_{x_n \to x_{-}} h_r(x_n) = 0  \; \; \text{or} \; {\lim}_{x_n \to x_{+}}h_r(x_n) = 0, r = A, T\}.$$  
We require a regularity condition on $h_A$ and $h_T$ in that at each $x \in S^o$, $h_T$ and $h_A$ converge to $0$ at certain rates in left and right neighborhoods or stay away from 0. We assume the cardinality of $S^o$, $|S^o| = m \; (m \geq 1)$, and  $S^o = \{s_1, \cdots, s_m\}$ for some $0 \leq s_1 < s_2 < s_m \leq 1$. More specifically, we assume there exists a small 
$\delta > 0$ such that the intervals $(s_j -\delta, s_j + \delta)$ are disjoint (in case $s_1 = 0$ or $s_m = 1$, the intervals are obviously truncated to $[0, \delta]$ and $[1 - \delta, 1]$, respectively) and 
\begin{itemize}
\item[{(i)}] for each $s_j \in S^o$, in its either left or right neighborhood, the source density $h_A$ converges to 0 at a polynomial rate, i.e.,  
\begin{eqnarray}
\label{eq:gene-da}
\underline{c}((x - s_j)_{-})^{a^A_{j, L}-1}  \leq h_A(x) &\leq& \overline{c} ((x - s_j)_{-})^{a^A_{j, L} - 1},  \text{for} \; \; x \in [s_j - \delta, s_j], \nonumber\\
\underline{c}((x - s_j)_{+})^{a^A_{j, R}-1}  \leq h_A(x) &\leq& \overline{c} ((x - s_j)_{+})^{a^A_{j, R} - 1},  \text{for} \; \; x \in [s_j,  s_j + \delta], 
\end{eqnarray}
with $\min(a^A_{j, L}, a^A_{j, R}) \geq 1$ and $\max(a^A_{j, L}, a^A_{j, R}) > 1$, which means that $h_A$ approaches 0 at least on one side of each singularity point in $S^o$.
For each $s_j$, the target density $h_T$ either converges to 0 at a polynomial rate in its left or right neighborhood, or stays away from $0$ on both sides, i.e.,  
\begin{eqnarray}
\label{eq:gene-dT}
\underline{c}((x - s_j)_{-})^{a^T_{j, L} - 1}  \leq h_T(x) &\leq& \overline{c} ((x - s_j)_{-})^{a^T_{j, L} - 1},  \text{for} \; \; x \in [s_j - \delta, s_j], \nonumber\\
\underline{c}((x - s_j)_{+})^{a^T_{j, R} - 1}  \leq h_T(x) &\leq& \overline{c} ((x - s_j)_{+})^{a^T_{j, R} - 1},  \text{for} \; \; x \in [s_j,  s_j + \delta], 
\end{eqnarray}
with $a^T_{j, L}, a^T_{j, R} \geq 1$.

\item[{(ii)}] Define $\mcB = \cup_{j=1}^m[s_j - \delta, s_j + \delta]$. For any point $x \in [0, 1]\backslash \mcB^c$, both $h_A$ and $h_T$ stay away from $0$, i.e.,  
\begin{eqnarray*}
\underline{c} \leq h_A(x), h_T(x) &\leq& \overline{c}.
\end{eqnarray*}
\end{itemize}

Define the set of transfer singularity points (TSP) as 
\[
S_{TSP} = \left\{s_j \in S^o: a^A_{j, L} > 1 + a^T_{j, L} (1 + \frac{1}{2\beta}) \; \text{or} \; a^A_{j, R} > 1 + a^T_{j, R} (1 + \frac{1}{2\beta}) \right\}.
\]

Let $J = \{(a^A_{j, r}, a^T_{j, r}): r = L, R, j = 1, \cdots, m, s_{j} \in S_{TSP}, a^A_{j, r} > 1 + a^T_{j, r}(1 + \frac{1}{2\beta})\}$. Clearly, the cardinality of 
$J$, $|J|$, is at least as the cardinality of $S_{TSP}$, $|S_{TSP}|$, and smaller than $2|S_{TSP}|$. Noting that if $s_1 = 0$ (or $s_m = 1$), 
then there is no left (or right) neighborhood, so the associated order parameters $(a^A_{1, L}, a^T_{1, L})$ (or $(a^A_{m, R}, a^T_{m, R})$) do not exist and must be excluded from $J$.  Let $\{x: a \leq x \leq b\} = \emptyset$ if $a > b$. When there is no TSP, then set $S_{TSP} = \emptyset$; otherwise, denote the pair $(a^A_{(1)}, a^T_{(1)}) := \arg\min_{((a^A_{j, r}, a^T_{j, r}) \in J)} \frac{2\beta + a^T_{j, r}}{2\beta + a^A_{j, r}} $. Theorem~\ref{thm:gen} below presents the minimax learning rate.

 \begin{theorem}
\label{thm:gen}
For $h_r, r = T, A$ satisfying the above two conditions (i) and (ii), we have that 
\begin{eqnarray*}
R_{h_T}(\mcH(\beta, K); \overline{\sigma}; \mcD, \mcD^T) \asymp \left\{
\begin{array}{ccc}
(n + n_T)^{-\frac{2\beta}{2\beta + 1}} &\forall n, n_T& S_{TSP} = \emptyset \\
n^{-\frac{2\beta + a^T_{(1)}}{2\beta + a^A_{(1)} }} & n_T \leq c_1 n^{\frac{2\beta + a^T_{(1)}}{2\beta + a^A_{(1)}}} & |S_{TSP}| \geq 1\\
n_T^{-\frac{2\beta}{2\beta + 1}} \left( \frac{n_T}{n} \right)^{\frac{1 + (a^T_{(1)} - 1)/(2\beta + 1)}{a^A_{(1)} - 
a^T_{(1)}}} &c_2 n^{\frac{2\beta + a^T_{(1)}}{2\beta + a^A_{(1)}}} \leq n_T \leq c_3 n &   |S_{TSP}| \geq 1 \\
n_T^{-\frac{2\beta}{2\beta + 1}}& c_4 n \leq n_T &  |S_{TSP}| \geq 1.\\
\end{array}
\right.
\end{eqnarray*}
\end{theorem}

\textbf{Proof of Theorem~\ref{thm:gen}}. We show the detailed proof procedure on the minimax lower bound next and leave the proof of upper bound in the SM Section S.1.1. 

We first focus on the case $h(x) \sim Beta(a^A, 1)$ and 
$h_T(x) \sim Beta(a^T, 1)$.  The extension to the general case can be seen readily afterwards. Similar to the Proof of 
Theorem~\ref{lem:ratef-beta-lower}, we establish the lower bound on case (I): $a^A \leq 1 + a^T(1 + \frac{1}{2\beta})$ and case (II)
$a^A > 1+ a^T(1 + \frac{1}{2\beta})$, separately. 

Case (I): $a^A \leq 1 + a^T(1 + \frac{1}{2\beta})$. Similarly to before, the lower rate $(n + n_T)^{-\frac{2\beta}{2\beta + 1}}$ again holds generally regardless of $\beta$, 
$n$ and $n_T$.

Case (II): $a^A > 1+ a^T(1 + \frac{1}{2\beta})$. As in the proof of Theorem~\ref{lem:ratef-beta-lower}, we only need to handle situations 1 and 2, as done below.  With 
$h_T \sim Beta(a^T, 1)$ replacing the uniform distribution, we have 
\[
\int_0^1 \psi^2_j(x)a^T x^{a^T - 1} dx \leq \frac{a_T A^2}{m^{2\beta + a^T}}j^{a^T - 1}. 
\]
From the proof of Theorem~\ref{lem:ratef-beta-lower}, for $f_t(x) = \sum_{j=1}^{J}t(j)\psi_j(x)$, we have 
\[
D(P^T_{f_t} || P_{f_0}^T) \leq \frac{1}{2\sigma^2} \frac{a^T A^2}{m^{2\beta + a^T}} J^{a^T}.
\]
Consequently, we know $f_t \in \mcF_{\eta}$ for all $t \in \mcT_J$ if 
$\eta^2 \geq \frac{c}{2\sigma^2}\left( n \frac{J^{a^A}}{m^{2\beta + a^A}} + n_T \frac{J^{a^T}}{m^{2\beta + a^T}}\right)$ for some constant $c > 0$. 

Similarly as before, in $\mcT_J$, by Gilbert-Varshamov bound, there exists at least $2^{c' J}$ elements such that for any pair, we have 
\[
\int_0^1 \left( f_{t_1}(x) - f_{t_2}(x) \right)^2a^T x^{a^T - 1} dx \geq c'' \frac{J^{a^T}}{m^{2\beta + a^T}}, 
\]
where $c'$ and $c''$ are constants that do not depend on $n$ or $n_T$.

Situation 1: $n_T \leq c_1 n^{\frac{2\beta + a^T}{2\beta + a^A}}$ for some $c_1 > 0$. By taking $J$ as a large enough constant and 
$m \asymp n^{\frac{1}{2\beta + a^A}}$, as in the proof of Theorem~\ref{lem:ratef-beta-lower}, the minimax lower bound holds. 

Situation 2: $c_1 n^{\frac{2\beta + a^T}{2\beta + a^A}} \leq n_T \leq c_2 n$ for some constants $0 < c_1, c_2 < \infty$. In this case, we take 
$m \asymp n^{\frac{1}{2\beta + 1}}(n_T/n)^{\frac{a^A - 1}{(a^A - a^T)(2\beta + 1)}}$ and 
$J \asymp m (n_T/n)^{\frac{1}{a^A - a^T}}$. Then correspondingly, we take 
\[
\eta^2 \asymp \frac{n m^{a^A}(n_T/n)^{\frac{a^A}{a^A - a^T}} }{m^{2\beta + a^A}},
\]
and the associated lower bound from Fano's inequality is of order 
\[
\frac{J^{a^T}}{m^{2\beta + a^T}} \asymp n_T^{-\frac{2\beta}{2\beta + 1}}(n_T/n)^{\frac{1 + \frac{a^T - 1}{2\beta + 1}}{a^A - a^T}},
\]
as can be directly verified. This completes the proof of the minimax lower rates in the case $h_T \sim Beta(a^T, 1)$ and 
$h_A \sim Beta(a^A, 1)$. 

Next, we extend the result to the more general setting where we have a single potential singularity point at $0$. More specifically, assume there 
exist $0 < \delta < 1$ and $ 0 < \underline{c} < \overline{c} < \infty$ such that 
$\underline{c} x^{a^A - 1} \leq h_A(x) \leq \overline{c}x^{a^A - 1}$, 
$\underline{c}x^{a^T - 1} \leq h_T(x) \leq \overline{c}x^{a^T - 1}$ for some $0 \leq x \leq \delta$ and 
$\inf_{\delta \leq x \leq 1} h_A(x)$ and $\inf_{\delta \leq x \leq 1} h_T(x)$ are positive. Here $\max(a^A, a^T) > 1$.  With this setup, it is straightforward to derive that for $J/m \to 0$, we have 
\[
D(P^T_{f_t} \| P^T_{f_0}) \leq \frac{c}{2\sigma^2}\frac{J^{a^T}}{m^{2\beta + a^T}}, \; 
D(P_{f_t} \| P_{f_0}) \leq \frac{c}{2\sigma^2}\frac{J^{a^A}}{m^{2\beta + a^A}},
\]
and as before, in $\mcT_{J}$, there exists at least $2^{cJ}$ elements such that for any pair, 
\[
\int_0^1 \left( f_{t_1}(x) - f_{t_2}(x) \right)^2 h_T(x) dx \geq \tilde{c} \frac{J^{a^T}}{m^{2\beta + a^T}},
\]
for some positive constants $c$ and $\tilde{c}$. The rest of the proof goes through similarly as before. 

Next, suppose we have a general potential singularity point $0 < s \leq 1$ on the left. Then the functions 
$\psi_j$ are modified to be:
\[
\psi_j(x) = m^{-\beta}\psi\left( m(s-x) - (j-1) \right),
\] 
for $1 \leq j \leq J$. The rest of the proof remains the same. If the potential singularity is on the right, the derivation can be done the same way. 
Finally, if we have multiple potential singularity points on the left or right, we apply the lower bounding argument for 
each (left or right) individually. Then the worst of the lower bounds becomes the overall minimax lower bound. 

Based on the results from Theorem~\ref{thm:gen}, we can see that when 
$S_{TSP} = \emptyset$, the minimax learning rate reduces to the classical optimal rate, given by $(n + n_T)^{-\frac{2\beta}{2\beta + 1}}$. When there is exactly one TSP, the minimax learning rate depends on the relationship between $n$ and $n_T$, which in turn is determined by the order parameters $a^T_{1, r}, a^A_{1, r}, r = L, R$ and the smoothness parameter $\beta$. In this case, the SLP happens if and only if $n_T \in  \left[c_1n^{\frac{2\beta + a^T_{(1)}}{2\beta + a^A_{(1) }}}, c_2n\right]$.  When multiple TSPs are present, the overall minimax learning rate is determined by the worst learning rate associated with any of the TSPs. Thus, if the conditions for SLP are met at least one TSP, then SLP occurs, but the degree of improvement depends on the worst learning rate among all TSPs. Furthermore, the improved rate increases as 
$a^A_{(1)}$ decreases and $a^T_{(1)}$ increases. 

\subsection{An illustration}
Next, we consider a simple case with two transfer singularity points including one left end point ($s_1 = 0$) and one internal point ($s_2 = 0.5$) to further illustrate the results of Theorem~\ref{thm:gen}. An example class of $(h_A, h_T)$ that may lead to the two TSPs is as follows:
\begin{eqnarray*}
h_r(x) &=&\left\{ \begin{array}{cc}
c_{r, 1}x^{a^r_{1,R}-1}(0.5 - x)^{a^r_{2, L}-1} & 0 \leq x \leq 0.5 \\
c_{r, 2}(x - 0.5)^{a^r_{2, R}-1} & 0.5 \leq x \leq 1,\\
 \end{array}
 \right. \; \; \; \text{for} \; \; r = A, T,
\end{eqnarray*}
where $c_{r, 1}, c_{r, 2}$ are positive constants selected such that the normalization condition $\int_{0}^1 h_r(x) dx = 1$ holds for both $r = A$ and $r = T$.
The parameters $\{a^r_{1, R}, a^r_{2, L}, a^r_{2, R}, r = A, T\}$ are all positive constants that are no less than 1. These parameters are additionally required to satisfy  
$a^A_{1, R} > 1 + a^T_{1, R}(1 + \frac{1}{2\beta})$, as well as the condition that $a^A_{2, r} > 1 + a^T_{2, r}(1 + \frac{1}{2\beta})$ holds for at least one value of $r \in \{L, R\}$. Under these conditions, $x = 0$ and $x = 0.5$ constitute the two TSPs. Particularly,  we consider the following two special cases, where case (1) represents that an extra TSP at $x = 0.5$ does not affect the improved rate, and case (2) represents that an extra TSP at $x = 0.5$ damages the improved rate.

(1) $\beta = 0.5, a^T_{1, R} = 1.5, a^A_{1, R} = 4.5, a^T_{2, R} = 2, a^A_{2, R} = 5.5, a^T_{2, L} = 2, a^A_{2, L} = 2$. 
We can obtain that  $J = \{(4.5, 1.5), (5.5, 2)\}$ and $(a^A_{(1)}, a^T_{(1)}) = (4.5, 1.5)$. Then, it follows that 
\begin{eqnarray*}
R_{h_T}(\mcH(\beta, K); \overline{\sigma}; \mcD, \mcD^T) \preceq \left\{
\begin{array}{cc}
n^{-2.5/5.5}& n_T \preceq n^{\frac{2.5}{5.5}}\\
n_T^{-\frac{1}{2}} \left( \frac{n_T}{n} \right)^{\frac{1.25}{3}}&n^{\frac{2.5}{5.5}} \preceq n_T \preceq n \\
n_T^{-1/2}& n \preceq n_T.\\
\end{array}
\right.
\end{eqnarray*}
This learning rate is the same as those with only one TSP at $x = 0$. Thus, the extra TSP does not affect the improved rate when SLP happens. 

(2) $\beta = 0.5, a^T_{1, R} = 1.5, a^A_{1, R} = 4.5, a^T_{2, R} = 1.5, a^A_{2, R} = 5, a^T_{2, L} = 2, a^A_{2, L} = 2$. We can obtain that 
$J = \{(4.5, 1.5),  (5, 1.5)\}$ and $(a^A_{(1)}, a^T_{(1)}) = (5, 1.5)$. Then, it follows that 
\begin{eqnarray*}
R_{h_T}(\mcH(\beta, K); \overline{\sigma}; \mcD, \mcD^T) \preceq \left\{
\begin{array}{cc}
n^{-2.5/6}& n_T \preceq n^{\frac{2.5}{6}}\\
n_T^{-\frac{1}{2}} \left( \frac{n_T}{n} \right)^{\frac{1.25}{3.5}}&  n^{\frac{2.5}{6.5}} \preceq n_T \preceq n \\
n_T^{-1/2}& n \preceq n_T.\\
\end{array}
\right.
\end{eqnarray*}

\begin{figure}
\centering
\includegraphics[height=6cm, width=13cm]{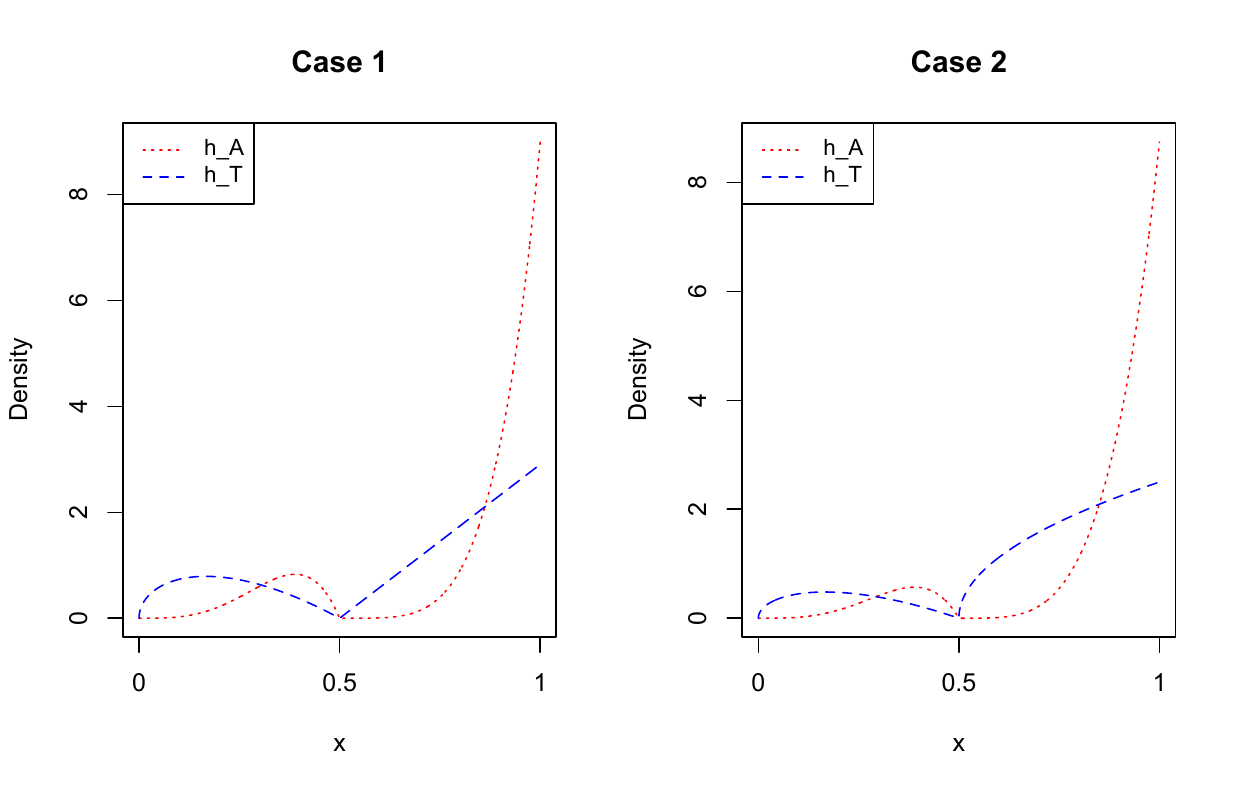}
\caption{One example of the two density functions $(h_A, h_T)$ under two cases, where case (1): $h_A(x) = \frac{1}{0.005}\left( x^{3.5}(0.5 - x)I_{[0, 0.5]}(x) + (x - 0.5)^{4.5}I_{[0.5, 1]}(x)\right)$ (dotted red line), $h_T(x) = \frac{1}{0.172}\left(x^{0.5}(0.5-x)I_{[0, 0.5]}(x) + (x - 0.5)I_{[0.5, 1]}(x)\right)$ (dashed blue line); and case (2): $h_A(x) = \frac{1}{0.007}\left(x^{3.5}(0.5 - x)I_{[0, 0.5]}(x) + (x - 0.5)^4I_{[0.5, 1]}(x)\right)$ (dotted red line), $h_T(x) = \frac{1}{0.283}\left(x^{0.5}(0.5 - x)I_{[0, 0.5]}(x) + (x - 0.5)^{0.5}I_{0.5, 1}(x)\right)$ (dashed blue line).}
\label{fig-ill}
\end{figure}

In this case, the learning rate with the extra TSP at 0.5 is changed and SLP may still happen but with a less improved rate. Figure~\ref{fig-ill} displays plots of the two density functions ($h_A, h_T$) under the above two cases with $(c_{A, 1}, c_{A, 2}, c_{T, 1}, c_{T,2}) = (1/0.005, 1/0.005, 1/0.172, 1/0.172)$ in case (1) and $(c_{A, 1}, c_{A, 2}, c_{T,1}, c_{T, 2}) = (1/0.007, 1/0.007, 1/0.283, 1/0.283)$ in case (2). 


\subsection{Densities with unknown parameters}

In practice, when we do not know the true $h_T$ and $h_A$, 
it is worthy of noting that consistent estimates of $h_T$ and $h_A$ may not be enough to ensure the same rate of convergence when knowledge of $h_A$ and $h_T$ is used to construct the regression estimator. However, if the density function is known up to a finite number of parameters, the minimax learning rate may still be achieved. 

More specifically, let's consider $h_T(x) \sim Beta(a_1^T, a_2^T)$ and $h_A(x) \sim Beta(a_1^A, a_2^A)$ where $a^T_1, a^T_2, a^A_1, a^A_2>0$ are unknown parameters. 
We next show that using maximum likelihood estimation (MLE) to estimate $a^T_1, a^T_2, a^A_1, a^A_2$ does not affect the TLR. For simplicity in derivation, we consider dividing both the target and source data into two halves. The first half will be used to estimate $a^T_1, a^T_2$ and $a^A_1, a^A_2$, respectively, while the second halves will be used to estimate $f$. We denote the density estimators based on MLE from the respective first half of the data as 
$\hat{h}_T := h_T(\cdot; \hat{\btheta}_T)$ and $\hat{h}_A := h_A(\cdot; \hat{\btheta}_A)$, where $\btheta_r = (a^r_1, a^r_2)^\top$. 
Consider a compact parameter space $\Theta := \{ (a^A_1, a^A_2, a^T_1, a^T_2): \underline{a} \leq a^r_j \leq \overline{a}, r = A, T, j = 1, 2\}$ for some known constants 
$0 < \underline{a} < \overline{a} < \infty$.  Let $\hat{\btheta} := (\hat{\btheta}_A, \hat{\btheta}_T)$ be estimators obtained based on MLE from the first half of the data, and 
 $ \hat{f}_{n, n_T}(\cdot; \hat{\btheta})$ be the estimators (NW estimator for $\beta \in (0, 1]$ and local polynomial estimator
 for $\beta > 1$, considered in Theorem~\ref{lem:ratef-beta} but now the bandwidth obtained using the estimated $\hat{h}_T$ and $\hat{h}_A$), based on the second halves of the target and source data.
Particularly, write $(x - t_n(x))^* = \max(x - t_n(x), 0)$ and $(x + t_n(x))^* = \min(x + t_n(x), 1)$.
Let $\hat{H}_r(x \pm t_n(x)) = \int_{(x-t_n(x))^*}^{(x + t_n(x))^*} \hat{h}_r(s)ds$ for $r = A, T$.  Let $\tilde{H}\left( x \pm t_n(x) \right) = \left(n H_A\left( x \pm t_n(x) \right) + 
n_T H_T(x \pm t_n(x))\right)/(n + n_T)$ and $\hat{\tilde{H}}\left( x \pm t_n(x) \right) = \left(n \hat{H}_A\left( x \pm t_n(x) \right) + 
n_T \hat{H}_T(x \pm t_n(x))\right)/(n + n_T)$. 

Define 
$\hat{\mcE}_n = \left\{\lambda_{\min}\left( \sum_{i=1}^{n + n_T}M_i \right)
\geq (n+ n _T) \hat{\tilde{H}}\left( x \pm t_n(x) \right) \frac{1}{16T_0} \right\},$
 where $T_0$ is some constant depending on $\overline{a}, \underline{a}$, 
$\lambda_{\min}(A)$ represents the minimum eigenvalue of $A$, $M_i = \mcZ\left( \frac{X_i - x}{t_n(x)} \right)\times $ $\mcZ^\top\left( \frac{X_i - x}{t_n(x)} \right)I(|X_i - x| < t_n(x))$, and 
 $\mcZ(u) = \left(1, u, u^2/2!, \cdots, u^{l}/l! \right)^\top$. Let 
 \[\tilde{f}_{n, n_T}(x; \hat{\btheta}) = I(\hat{\mcE}_n)\sum_{i=1}^{n + n_T} Y_i W^*_{ni}(x) \bm{I}\left( \sum_{i=1}^{n + n_T}I(|X_i - x| \leq t_n(x)) > 0 \right), \]
 where 
 \begin{eqnarray*}
 W^*_{ni}(x) &=& \frac{1}{(n + n_T) t_n(x)}\mcZ^\top(0) \mcB^{-1}_{nx} \mcZ\left( \frac{X_i - x}{t_n(x)} \right)I \left( |X_i - x| \leq t_n(x) \right), \\
 \mcB_{nx} &=& \frac{1}{( n + n_T) t_n(x)}\sum_{i=1}^{ n + n_T} M_i.
 \end{eqnarray*}
 Then, the final estimator takes the following form
 \begin{eqnarray}
 \label{eq:de}
 \hat{f}_{n, n_T}(x; \hat{\btheta}) &=& \tilde{f}_{n, n_T}(x; \hat{\btheta})I\left( \tilde{f}_{n, n_T}(x; \hat{\btheta}) < T_1 \right), 
 \end{eqnarray}
where $T_1$ is some large constant, depending on $\overline{\sigma}, \beta, K$. In practice, the value $t_n(x)$ can be calculated based on formula~\eqref{eq:lb} with $H_T$ and $H_A$ replaced by 
 $\hat{H}_T$ and $\hat{H}_A$, and the constant $T_1$ can be a sufficiently large constant, and $T_0 = 2 T_1$.


 \begin{theorem}
\label{thm:est}
Assume $h_T(x) \sim Beta(a_1^T, a_2^T)$ and $h_A(x) \sim Beta(a_1^A, a_2^A)$ where $(a^A_1, a^A_2, a^T_1, a^T_2) := \btheta \in \Theta$. 
Then, the aforementioned estimator satisfies
\[
\sup_{f \in \mcH(\beta, \kappa), 0 < \sigma \leq \overline{\sigma}} R_{h_T}(f; \hat{f}_{n, n_T}(\cdot; \hat{\btheta}); n, n_T)= O\left(R_{h_T}(\mcH(\beta, K); \overline{\sigma}; \mcD, \mcD^T)\right).
\]
\end{theorem}

The proof of Theorem~\ref{thm:est} is provided in SM Section S.1.2.
\subsection{Adaptation to the smoothness parameter $\beta$}

The above mentioned estimator $\hat{f}$ is not adaptive to the smoothness parameter $\beta$. An adaptive version of this estimator can be obtained at the price of an extra logarithmic factor.
In particular, using the same strategy as illustrated in section 4 of \cite{qian2016randomized}, which is constructed following the idea in \cite{lepskii1991problem}, we can obtain a proper estimator $\hat{\beta}$. We use the first half of the source data to estimate $\hat{\beta}$ when $n > n_T$; otherwise, we use the first half of  the target data for estimation. For simplicity, we assume $h_T(x) \sim Beta(a_1^T, a_2^T)$ and $h(x) \sim Beta(a_1, a_2)$ with $a_1^T, a_2^T, a_1, a_2$ unknown.  To ensure that $\hat{\beta}$ retains desirable estimation properties \cite[Proposition 4.1,][]{qian2016randomized}, 
it is desirable that the density function is bounded away from zero and infinity.  To achieve this, we use data with $x$ in $[\underline{c}, 1- \underline{c}]$, where $\underline{c}$ is a very small positive constant. Specifically, we construct $\hat{\beta}$ using $\mcD_c$ for the source data (or $\mcD_c^T$ for the target data), where $\mcD_c := \{(X_i, Y_i):  X_i \in [\underline{c}, 1 - \underline{c}], i = 1, \cdots, \lfloor n/2\rfloor\}$ and $\mcD^T_c := \{(X^T_i, Y^T_i): X^T_i \in [\underline{c}, 1 - \underline{c}], i = 1, \cdots, \lfloor n_T/2\rfloor\}$. Then the adaptive estimator, defined as $\hat{f}_{n, n_T}(x; \hat{\btheta}, \hat{\beta})$,  is the estimator in equation~\eqref{eq:de}, with $l$ replaced by $\lfloor \hat{\beta} \rfloor$. 

Let's assume that the source data is used to construct the estimate. For $\beta \in (\underline{\beta}, \overline{\beta})$, following \cite{qian2016randomized}, define two integers $\tau^* = \max\{\tau + 1: 2^{\tau} \leq n^{\frac{1}{2\underline{\beta} + 1}}\}$ and 
$\tau_* = \max\{\tau: 2^{\tau} \leq n^{\frac{1}{2\overline{\beta} + 1}}\}$. For any integer $\tau$, let 
$\mu_{\tau} = 2^{-\tau}$, and $\beta_{\tau}$ be the real number that satisfies $\mu_{\tau} = n^{-\frac{1}{2\beta_{\tau} + 1}}$. Clearly, 
there exists a constant $c > 0$ such that 
$\beta_{\tau} - \beta_{\tau + 1} \leq \frac{c}{\log n}$ for any 
$\tau \in [\tau_*, \tau^*]$. 
Given $\tau$, 
we evenly partition the domain into $1/\mu_{\tau}$ bins with bin width $(1 - 2\underline{c})\mu_{\tau}$, and let 
$\mcD_{\tau}(x)$ denote the bin that contains $x \in [\underline{c}, 1-\underline{c}]$. Then, with any given $x \in [\underline{c}, 1-\underline{c}]$ and 
$\tau \in [\tau_*, \tau^*]$, we estimate $f(x)$ by $\hat{m}_{\tau}(x) = \frac{\sum_{ i \in H_{\tau}(x)} Y_i}{M_{\tau}(x)}$, where 
$H_{\tau}(x) = \{i: X_i \in \mcD_{\tau}(x), (X_i, Y_i) \in \mcD_c\}$ when the source data is used to estimate $\hat{\beta}$, and $M_{\tau}(x)$ is the size of $H_{\tau}(x)$. Define 
$\hat{\tau}$ to be $
\min\{\tau \in [\tau_*, \tau^*]: \|\hat{m}_{\tau} - \hat{m}_{\tau_2}\|_{\infty} \leq c \mu^{\beta_{\tau_2}}_{\tau_2}\log n \;
\text{for every} \; \tau_2 \; \text{satisfying} \; \tau < \tau_2 \leq \tau^*
\},
$ 
where $c$ is a constant satisfying $c > 4\kappa$. Then $\hat{\beta}$ can be estimated by 
$\hat{\beta} = \min(\beta_{\hat{\tau}} - \frac{c \log \log n}{\log n}, \overline{\beta})$, where 
$c$ is a constant satisfying $c > \frac{(2\overline{\beta} + 1)^2}{2 \underline{\beta}}$.

Consider a sub-class $\mcH_{0, \underline{c}}(\beta, \kappa)$ of $\mcH(\beta, \kappa)$ as follows. Given $\tau \in \mbN$ and 
$x \in [\underline{c},1-\underline{c}]$, define 
$K_{\tau, h} f(x) =: \mbE\left[ f(X)\mid X \in \mcD_{\tau}(x) \right] = \int_{D_{\tau}(x)}f(t)h(t)dt/ \int_{D_{\tau}(x)} h(t)dt.$
Then 
\begin{eqnarray*}
&&\mcH_{0, \underline{c}}(\beta, \kappa) =: \left\{ 
f \in \mcH(\beta, \kappa): \; \text{there exists} \; 0 < \kappa_1 < \kappa \; \text{and} \; \tau_0 > 0 \; \text{such that for every} \; \tau \geq \tau_0 \right.\\
&&\hspace{0.5cm} \left.
\min\left(\sup_{x \in [\underline{c},1-\underline{c}]}|K_{\tau, h}f(x) - f(x)|, \sup_{x \in [\underline{c},1-\underline{c}]}|K_{\tau, h_T}f(x) - f(x)|\right) > \kappa_1 n^{-\frac{\beta}{2\beta_{\tau} + 1}} \; 
\right\}.
\end{eqnarray*}
It can be seen that for any $f \in \mcH_{0, \underline{c}}(\beta, \kappa) $, $f \not\in \mcH(\tilde{\beta}, \kappa)$ for every 
$\tilde{\beta} > \beta$. More importantly, from Proposition 4 in \cite{gine2010confidence}, we can see that 
$\mcH(\beta, \kappa) \backslash \mcH_{0, \underline{c}}(\beta, \kappa)$ is \emph{nowhere dense} in the relevant H$\ddot{\text{o}}$lder-norm topologies.

\begin{theorem}
\label{thm:ad}
Under the conditions of Theorem~\ref{thm:est}, if $\beta \in (\underline{\beta}, \overline{\beta})$ for two positive constants $\underline{\beta}, \overline{\beta}$ satisfying $0 < \underline{\beta} < \overline{\beta} < \infty$, then the estimator  $\hat{f}_{n, n_T}(x; \hat{\btheta}, \hat{\beta})$ satisfies 
\[
\sup_{\substack{f \in \mcH_{0, \underline{c}}(\beta, \kappa) \\ 0 < \sigma \leq \overline{\sigma}}} R_{h_T}(f; \hat{f}_{n, n_T}(\cdot; \hat{\btheta}; \hat{\beta}); n, n_T)= O\left(R_{h_T}(\mcH(\beta, K); \overline{\sigma}; \mcD, \mcD^T) \big(\log (n + n_T)\big)^c\right),
\]
where $c$ is a constant depending on $\underline{\beta}, \overline{\beta}$. 
\end{theorem}

The proof of Theorem~\ref{thm:ad} is provided in the SM Section S.1.3.

\section{Concluding Remarks}
\label{sec:con}

In this work, we find a new synergistic learning phenomenon that under the framework of domain adaptation,  when the ratio of $h_T(x)$ and 
$h(x)$ approaches $\infty$ sufficiently fast, the minimax convergence rate based on the source data and the target data together becomes faster than the better rate of convergence based on the source data or target data only, as suggested by the fascinating upper bounds obtained by \cite{schmidt2022local}.
And then, the synergistic learning phenomenon occurs when and only when the target sample size is significantly smaller than but not too much smaller than the source sample size. The SLP materializes in two different ways: One is that the target data help 
alleviate the difficulty in estimating the regression function at points where $h(x)$ is close to zero and the other is that the source 
data help the estimation at majority of the points where $h(x)$ is not very small due to its larger size than $n_T$. In real application of domain adaptation, it is perhaps expected that the target and source density ratio can be very large in various local regions. The theoretical understanding obtained in this work that allows multiple singularity points sheds some light on how the source and target data can be used wisely to create a synergy. 

For simplicity, this work focuses on regression functions in the one-dimensional H$\ddot{\text{o}}$lder class. An important direction for future research is understanding the behavior when densities and regression functions belong to general multi-dimensional families, particularly when $h_T$ and $h$ require nonparametric estimation.


\appendix 
\renewcommand{\theequation}{A.\arabic{equation}}
\renewcommand{\thesection}{A.\arabic{section}}
\renewcommand{\thesubsection}{A.\arabic{section}.\arabic{subsection}}
\renewcommand{\thetheorem}{A.\arabic{theorem}}
\renewcommand{\thelemma}{A.\arabic{lemma}}
\renewcommand{\thecorollary}{A.\arabic{corollary}}
\renewcommand{\theexample}{A.\arabic{example}}
\renewcommand{\thefigure}{A.\arabic{figure}}
\renewcommand{\theproposition}{A.\arabic{proposition}}

\section{Existence of $t_n(x)$ and its properties}
\begin{lemma}
\label{lem:tn}
Let $h_T(x) \sim \text{Unif}(0, 1)$ and $h(x) = a x^{a-1}I( 0 \leq x \leq 1)$ for some $a > 0$. Assume $\beta > 0$. For any $n + n_T > 1$ and any $x \in [0, 1]$, there exist a unique solution $t_n(x)$ of the equation 
\[
t_n(x)^{2\beta} \left( \frac{n}{n + n_T} H(x \pm t_n(x))  + \frac{n_T}{n + n_T} H_T(x \pm t_n(x))\right) = \frac{1}{n + n_T}.
\]
Therefore the function $x \mapsto t_n(x)$ is well defined on $[0, 1]$. From now on, we refer to $t_n$ as the spread function. 
\end{lemma}

\begin{lemma}
\label{lem:tn1}
We have the following results:
\begin{itemize}
\item[{(i)}] $x \mapsto t_n(x)$ is $\mcH(1, 1)$.
\item[{(ii)}] for $n + n_T \geq 2^{2\beta}$, there exist unique solutions $x_1$ and $x_2$ to, respectively,
$t_n(x) = x$ and $t_n(x) = 1 - x$ satisfying $0 < x_1 < 1/2 < x_2 < 1$.
\item[{(iii)}] $t_n$ is differentiable on $(0, 1) \backslash \{x_1, x_2\}$ and the derivative $\dot{t}_n(x)$ is given by 
\begin{eqnarray*}
\frac{ \tilde{h}(x - t_n(x))\bm{I}(t_n(x) < x) - \tilde{h}(x + t_n(x))\bm{I}(t_n(x) < 1 - x)}{2\beta/((n + n_T)t_n(x)^{2\beta + 1}) + \tilde{h}(x + t_n(x))\bm{I}(t_n(x) < 1 - x) + \tilde{h}(x - t_n(x))\bm{I}(t_n(x) < x)},
\end{eqnarray*}
where $\tilde{h}(x) = \frac{n}{n + n_T}h(x) + \frac{n_T}{n + n_T}h_T(x)$.
\end{itemize}
\end{lemma}

The proofs of Lemmas~\ref{lem:tn}-\ref{lem:tn1} are provided in SM Section S.2.  

\section{Inequalities for NW and local polynomial estimators}
For $\beta \in (0, 1]$, we consider the Nadaraya-Watson (NW) estimator to estimate $f$. Denote $U(x) = \sum_{i=1}^{n + n_T}\bm{I}\{X_i \in (x \pm t_n(x))\}$. In particular, the NW estimator $\hat{f}_{n, n_T}$ takes the following form:
\begin{eqnarray}
\label{eq:nw}
\hat{f}_{n, n_T}(x) =  \frac{\sum_{i=1}^{n + n_T} Y_i\bm{I}\{X_i \in (x \pm t_n(x))\} }{\sum_{i=1}^{n + n_T}\bm{I}\{X_i \in (x \pm t_n(x))\}}\bm{I}(U(x) > 0).
\end{eqnarray}
This Nadaraya-Watson estimator is a local constant least squares approximation of the output. 

For $f \in \mcH(\beta, \kappa)$ with $\beta > 1$, $l = \lfloor \beta \rfloor$, we consider the local polynomial estimator. In particular,
 let $\mcZ(u) = \left(1, u, \cdots, u^{l}/l! \right)^\top$ and $\theta(x) = \left( f(x), f^{(1)}(x)t_n(x), f^{(2)}(x)t^2_n(x), 
 \cdots, \right.$ $\left. f^{(l)}(x)t^l_n(x) \right)^{\top}$ with $f^{(k)}(\cdot)$ being the $k$th derivative of $f(\cdot)$.  
 The local polynomial estimator of order $l$ of $\theta(x)$ is defined by 
 \[
 \hat{\theta}_{n, n_T}(x) = \arg\min_{\theta \in R^{l + 1}} \sum_{i=1}^{n + n_T} \left[ Y_i - \btheta^\top \mcZ\left( \frac{X_i - x}{t_n(x)} \right) \right]^2 I\left( |X_i - x| \leq t_n(x) \right).
 \]
 Then, the local polynomial estimator of order $l$ of $f(x)$ is defined by $\hat{f}_{n, n_T}(x) = \mcZ^\top(0) \hat{\theta}_{n, n_T}(x)$. 
Let $M_i = \mcZ\left( \frac{X_i - x}{t_n(x)} \right) \mcZ^\top\left( \frac{X_i - x}{t_n(x)} \right)  I(|X_i - x| < t_n(x))$. After simply calculation, we have 
 \begin{eqnarray}
 \label{eq:lp}
\hat{f}_{n, n_T}(x) = \sum_{i=1}^{n + n_T} Y_i W^*_{ni}(x) \bm{I}\left( \sum_{i=1}^{n + n_T}I(|X_i - x| \leq t_n(x)) > 0 \right), 
 \end{eqnarray}
 with 
 \begin{eqnarray*}
 W^*_{ni}(x) &=& \frac{1}{(n + n_T) t_n(x)}\mcZ^\top(0) \mcB^{-1}_{nx} \mcZ\left( \frac{X_i - x}{t_n(x)} \right)I \left( |X_i - x| \leq t_n(x) \right), \\
 \mcB_{nx} &=& \frac{1}{( n + n_T) t_n(x)}\sum_{i=1}^{ n + n_T} M_i.
 \end{eqnarray*}

The next Lemma~\ref{lem:upb} and Lemma~\ref{local-p} show the upper bound for NW estimator and local polynomial estimator, respectively.

\begin{lemma}
\label{lem:upb}
For each $x \in [0, 1]$ almost surely, the Nadaraya-Watson estimator $\hat{f}_{n, n_T}$ in expression~\eqref{eq:nw} satisfies 
\[
\mbE\left[ \big(\hat{f}_{n, n_T}(x) - f(x)\big)^2 \right] \leq \kappa^2 \big(t_n(x)\big)^{2\beta} + \frac{2\sigma^2 + \|f\|^2_{\infty}}{nH(x \pm t_n(x)) + n_TH_T(x \pm t_n(x))},
\]
where the expectation is taking over the observations $\{(X_i, Y_i)\}^{n + n_T}_{i=1}$ and $X_i \sim h(\cdot)$ for $i = 1, \cdots, n$ and 
$X_i \sim h_T(\cdot)$ for $i = n+1, \cdots, n+n_T$.
\end{lemma}

For $\beta \geq 1$, let $\hat{f}_{n, n_T}$ be the local polynomial estimator in~\eqref{eq:lp} of order $l$ of $f$ with $l = \lfloor \beta \rfloor$. Denote $\mcE_n := \left\{X_1, \cdots, X_{n+n_T}: \lambda_{\min}\left( \sum_{i=1}^{n + n_T}M_i \right)
\geq (n+ n _T) \tilde{H}\left( x \pm t_n(x) \right)/8 \right\}$, and
$\tilde{H}\left( x \pm t_n(x) \right) = \left(n H\left( x \pm t_n(x) \right) + 
n_T H_T(x \pm t_n(x))\right)/(n + n_T),$ where 
$\lambda_{\min}(A)$ represents the minimum eigenvalue of $A$. Let 
 $$ 
\tilde{f}_{n, n_T}(x) = \hat{f}_{n, n_T}(x) I(\mcE_n).
 $$
 For notation simplicity, we still use $\hat{f}_{n, n_T}(x)$ to replace the truncated $\tilde{f}_{n, n_T}(x)$ unless it is needed.

 \begin{lemma}
 \label{local-p}
 For all $x \in [0, 1]$ and a constant $c$, the following inequality holds:
\[
\mbE \left[ (\hat{f}_{n, n_T}(x) - f(x))^2 \right] \leq ct^{2\beta}_n(x) + \frac{\|f\|^2_{\infty} + 8 \sigma^2 (l+1)^2}{nH(x \pm t_n(x)) 
+ n_TH_T(x \pm t_n(x))},
\]
where the expectation is taking over the observations $\{(X_i, Y_i)\}_{i=1}^{n + n_T}$, 
and $X_i \sim h(\cdot)$ for $i = 1, \cdots, n$,  
and $X_i \sim h_T(\cdot)$ for $i = n + 1, \cdots, n + n_T$.
  \end{lemma}
 
 The proofs of Lemmas~\ref{lem:upb}-\ref{local-p} are provided in SM Section S.3.

\bibliography{ref}     

\end{document}